\documentclass[english,plain]{article}

\RequirePackage{xcolor,graphicx}

\RequirePackage{amsmath,amsfonts,amssymb,amsthm,mathtools,accents}
\RequirePackage{bm}

\RequirePackage[a4paper]{geometry}
\RequirePackage{subfig}
\usepackage{multirow}

\newcommand{\tr}{{\operatorname{tr}}}

\DeclareMathOperator*{\dof}{\operatorname{dof}}

\newcommand{\norm}[1]{{\left\Vert {#1} \right\Vert}}
\newcommand{\abs}[1]{{\left\vert {#1} \right\vert}}
\newcommand{\interior}[1]{\accentset{\circ}{#1}}

\usepackage[normalem]{ulem}

\graphicspath{{fig/}{./}}

\title{Performances of the mixed virtual element method on complex grids for underground flow}

\author{Alessio Fumagalli$^*$  \and Anna
Scotti$^*$ \and Luca Formaggia\thanks{MOX - Dipartimento di Matematica ``F. Brioschi'', Politecnico di
Milano, via Bonardi 9, 20133 Milan, Italy.
\texttt{alessio.fumagalli@polimi.it}, \texttt{anna.scotti@polimi.it},\  \texttt{luca.formaggia@polimi.it}.}}

\begin{document}

\maketitle

\begin{abstract}
    The numerical simulation of physical processes in the underground frequently
    entails challenges related to the geometry and/or data. The former are
    mainly due to the shape of sedimentary layers and the presence of fractures
    and faults, while the latter are connected to the properties of the rock
    matrix which might vary abruptly in space. The development of approximation
    schemes has recently focused on the  overcoming of such difficulties with
    the objective of obtaining numerical schemes with  good approximation
    properties.  In this work we carry out a numerical study on the
    performances of the Mixed Virtual Element Method (MVEM) for the solution of
    a single-phase flow model in fractured porous media. This method is able to
    handle grid cells of polytopal type and treat hybrid dimensional problems.
    It has been proven to be robust with respect to the variation of the
    permeability field and of the shape of the elements. Our numerical
    experiments focus on two test cases that cover several of the aforementioned
    critical aspects.
\end{abstract}

\section{Introduction}

The numerical simulation of subsurface flows is of paramount importance in many
environmental and energy related applications such as the management of
groundwater resources, geothermal energy production, subsurface storage of
carbon dioxide. The physical processes are usually modeled, under suitable
assumptions, by Darcy's law and its generalization to multiphase flow.

In spite of the simplicity of the Darcy model the simulation of subsurface flow is often
a numerical challenge due to the strong heterogeneity of the coefficients,
porosity and permeability of the porous medium, and to the geometrical
complexity of the domains of interest. At the spatial scale of reservoirs, or
sedimentary basins, the porous medium has a layered structure due to the
deposition and erosion of sediments, and tectonic stresses can create, over
millions or years, deformations, folds, faults and fractures. In realistic cases
the construction of a computational grid that honours the geometry of layers and
a large number of fractures is not only a difficult task, but can also give poor
results in terms of quality, creating, for instance, very small or
badly shaped elements in the vicinity of the interfaces.

In the framework of
Finite Volume and Finite Elements methods one possibility is to consider
formulations that allow for coarse and regular grids cut by the interfaces in
arbitrary ways. The Embedded Discrete Fracture Model, for instance,
\cite{Li2008,Panfili2014,Fumagalli2015},
can represent permeable fracture that cut the background grid by adding
additional trasmissibilities in the matrix resulting from the Finite Volumes
discretization; on the other hand the eXtended Finite Element Method can be used
to generalize a classical FEM discretization allowing for discontinuities inside
an element of the grid, see for example \cite{DAngelo2011,Formaggia2012,DelPra2015a,Flemisch2016}
for the application of this technique to
Darcy's problem.

A promising alternative consists in the use of numerical methods that are
robust in the presence of more general grids, in particular polygonal/polyhedral
grids and that impose mild restriction on element shape: this is the case for
the Virtual Element Method (VEM), introduced in
\cite{BeiraodaVeiga2014a,Brezzi2014,BeiraoVeiga2016} and successfully applied
now to a variety of problems, including elliptic problems in mixed form which is
the case of the Darcy model considered in this work. See also
\cite{Benedetto2014,Benedetto2016,Fumagalli2016a,Fumagalli2017,Fumagalli2017a}. By avoiding the explicit
construction of basis function VEM can indeed handle very general grids, which
might be useful in the aforementioned cases where the heterogeneity of the
medium and the presence of internal interfaces pose constraints to grid
generation. In the context of porous media simulations, mixed methods, i.e. methods that consider both velocity and pressure as unknowns of the problem, are of particular interest since they provide a good approximation of pressure as well as an accurate (and conservative) velocity field. For these reasons, we focus our attention on the Mixed Virtual Element Method (MVEM).

The aim of this work is to consider practical grid generation
strategies to deal with such complex geometries and to test the performances of
the MVEM method on the different types of grid proposed. In
particular, we want to investigate the impact of grid type and element shape on
properties of the linear system such as sparsity and condition number, and
eventually compare the errors. To this aim we will consider two test cases from
the literature, in particular two layers from the well-known  10$^{\rm th}$ SPE
Comparative Solution Project (SPE10) dataset, described in~\cite{Christie2001},
characterized by a complex permeability field, and a test case for fractured
media taken from~\cite{Flemisch2016a}. We focus our attention on grid generation strategies that
can be applicable in realistic cases: if it is certainly true that MVEM can
handle general polytopal grid the construction of such grids is often a
difficult task. For this reason, in addition to classical Delaunay triangular
grids we consider the case of Voronoi grids, rectangular Cartesian grids cut by
fractures, and grids generated by agglomeration. This latter strategy can be
applied as a post-processing to all other grid types with the advantage of reducing
the number of unknowns. For the numerical
implementation of the test cases we have used the publicly available library PorePy~\cite{Keilegavlen2019a}.

The paper is structured as follows: in Section \ref{sec:governing_eq} we present the mathematical
model, i.e. the single phase Darcy model in the presence of fractures
approximated as codimension 1 interfaces. Section \ref{sec:weak_form} is devoted
to the weak formulation of the problem just introduced. Section \ref{sec:vem}
introduces the
numerical discretization by the Virtual Element method, while in Section \ref{sec:grid_generation} we
describe the grid generation strategies used in the paper. Section \ref{sec:numerical_results} presents
the numerical tests, and Section \ref{sec:conclusion} is devoted to conclusions.

\section{Governing equations}\label{sec:governing_eq}

We now introduce the mathematical models considered in this
work. The realistic modeling of subsurface flows requires a complex set of non-linear equations and constitutive laws, however one of the key
ingredient (upon a suitable linearisation) is the single-phase flow model for a
porous media based on Darcy's law and mass conservation. We are here studying this model, keeping in mind that it might be seen as a part of a more complex model. In addition, it is of our interest to consider also fractures in the porous media, and this calls for a more sophisticated approach.

As already mentioned, we set our study in a saturated porous medium represented by the domain
$\Omega \subset \mathbb{R}^2$. The boundary of $\Omega$, indicated with $\partial \Omega$,
is supposed regular enough (e.g. Lipschitz continuous). The boundary is divided into two
disjoint parts $\partial_u \Omega$ and $\partial_p \Omega$ such that
$\interior{\partial_u \Omega} \cap
\interior{\partial_p \Omega} =
\emptyset$ and $\overline{\partial_u \Omega} \cup \overline{\partial_p \Omega} =
\overline{\partial \Omega}$. These portions of the
boundary will be used to define boundary conditions.

\subsection{Single-phase flow in the bulk domain}\label{subsec:darcy}

We briefly recall the mathematical model of single-phase flow in porous media, referring to classical results in literature, see \cite{Bear1972}, for details. We are interested in the computation of the vector field Darcy
velocity $\bm{u}$ and scalar field pressure $p$, which are solutions of the following
problem
\begin{align}\label{eq:darcy}
    \begin{aligned}
        &\begin{aligned}
            &\bm{u} + K\nabla p = \bm{0} \\
            &\nabla \cdot \bm{u} = f
        \end{aligned}
        && \text{in } \Omega,\\
        &\bm{u} \cdot \bm{n}_\partial = \overline{u} && \text{on }
        \partial_u \Omega,\\
        & p = \overline{p} && \text{on } \partial_p \Omega.
    \end{aligned}
\end{align}
The parameter $K$ is the $2 \times 2$ permeability tensor, which is
symmetric and positive definite. For simplicity, the dynamic viscosity of the
fluid is included into $K$. The source or sink term is named $f$. Finally, $\bm{n}_\partial$ is
the outward unit normal on $\partial \Omega$, $\overline{u}$ and $\overline{p}$ given boundary data.

We recall that the permeability tensor, for real applications, may vary several order of
magnitude from region to region (i.e., grid cells) and can be discontinuous.

\subsection{Fracture flow}\label{subsec:fracture}

We are interested in the simulation of single-phase flow in porous media
in the presence of fractures. For simplicity we start with one single
fracture. The model we are considering is the result of a
model reduction procedure that approximates the fracture as a lower dimensional
object and derives new equations and coupling conditions for the Darcy velocity and pressure both in the fracture and surrounding porous
medium. More details on this subject  can
be found in the following, not exhaustive, list of works
\cite{Alboin2002,Martin2005,DAngelo2011,Formaggia2012,Sandve2012,Brenner2016a,Fumagalli2017a,Antonietti,Scotti2017,Boon2018,Chave2018,Stefansson2018,Nordbotten2018,Berre2019b}.

\begin{figure}[tbp]
    \centering
\begingroup%
  \makeatletter%
  \providecommand\color[2][]{%
    \errmessage{(Inkscape) Color is used for the text in Inkscape, but the package 'color.sty' is not loaded}%
    \renewcommand\color[2][]{}%
  }%
  \providecommand\transparent[1]{%
    \errmessage{(Inkscape) Transparency is used (non-zero) for the text in Inkscape, but the package 'transparent.sty' is not loaded}%
    \renewcommand\transparent[1]{}%
  }%
  \providecommand\rotatebox[2]{#2}%
  \newcommand*\fsize{\dimexpr\f@size pt\relax}%
  \newcommand*\lineheight[1]{\fontsize{\fsize}{#1\fsize}\selectfont}%
  \ifx\svgwidth\undefined%
    \setlength{\unitlength}{154.02376791bp}%
    \ifx\svgscale\undefined%
      \relax%
    \else%
      \setlength{\unitlength}{\unitlength * \real{\svgscale}}%
    \fi%
  \else%
    \setlength{\unitlength}{\svgwidth}%
  \fi%
  \global\let\svgwidth\undefined%
  \global\let\svgscale\undefined%
  \makeatother%
  \begin{picture}(1,0.56888735)%
    \lineheight{1}%
    \setlength\tabcolsep{0pt}%
    \put(0,0){\includegraphics[width=\unitlength,page=1]{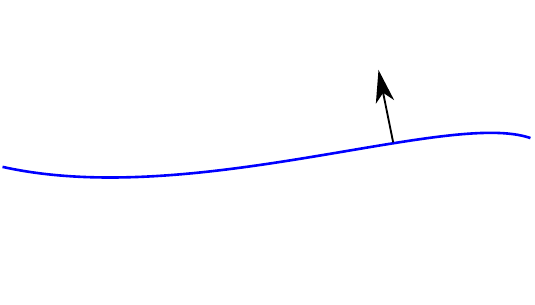}}%
    \put(0.73698117,0.3700154){\color[rgb]{0,0,0}\makebox(0,0)[lt]{\lineheight{1.25}\smash{\begin{tabular}[t]{l}$\bm{n}$\end{tabular}}}}%
    \put(0.38240354,0.18058272){\color[rgb]{0,0,0}\makebox(0,0)[lt]{\lineheight{1.25}\smash{\begin{tabular}[t]{l}$\gamma^-$\end{tabular}}}}%
    \put(0,0){\includegraphics[width=\unitlength,page=2]{fracture_reduced.pdf}}%
    \put(0.39773586,0.28375307){\color[rgb]{0,0,0}\makebox(0,0)[lt]{\lineheight{1.25}\smash{\begin{tabular}[t]{l}$\gamma^+$\end{tabular}}}}%
    \put(0.08275525,0.26558866){\color[rgb]{0,0,0}\makebox(0,0)[lt]{\lineheight{1.25}\smash{\begin{tabular}[t]{l}$\gamma$\end{tabular}}}}%
    \put(0.30096685,0.05899831){\color[rgb]{0,0,0}\makebox(0,0)[lt]{\lineheight{1.25}\smash{\begin{tabular}[t]{l}$\Omega^-$\end{tabular}}}}%
    \put(0.30096685,0.44758489){\color[rgb]{0,0,0}\makebox(0,0)[lt]{\lineheight{1.25}\smash{\begin{tabular}[t]{l}$\Omega^+$\end{tabular}}}}%
  \end{picture}%
\endgroup%

    \caption{Hybrid-dimensional representation of a fracture immersed in a porous
    media.}
    \label{fig:fracture_reduced}
\end{figure}

In the following the fracture is indicated with $\gamma$, and quantities related to
the porous media and the fracture are indicated with the subscript $\Omega$ and
$\gamma$, respectively. The fracture is described by a planar surface with normal
vector denoted by $\bm{n}$, which also defines a positive and negative side of
$\gamma$, indicated as $\gamma^+$ and $\gamma^-$, see Figure
\ref{fig:fracture_reduced} as an example.  Given a field $u$ in
$\Omega\setminus\gamma$ we indicate its trace on $\gamma^+$ and $\gamma^-$ as
$\tr u_+$ and $\tr u_-$, respectively.

The fracture is characterized by an aperture $\epsilon_\gamma$ which, in the
reduced model where the fracture has co-dimension one, is only a model parameter.
Finally, if the fracture touches the boundary we can apply natural or essential
given boundary conditions; we denote as $\partial_p \gamma$ and $\partial_u
\gamma$ the portions of $\partial \gamma$ where pressure and velocity are
imposed. We assume that $\interior{\partial_p \gamma} \cap \interior{\partial_u
\gamma} = \emptyset$ as well as ${\partial_p \gamma} \cup {\partial_u \gamma} =
{\partial \gamma}$.  If a fracture tip does not touch the physical boundary a
no-flow condition is imposed, so in this case we assume that the immersed tip
belongs to $\partial_u \gamma$ with an homogeneous condition.

We recall the system of equations that will be used
in the sequel. In the bulk porous medium $\Omega\setminus\gamma$ the problem is governed by the classic Darcy's equations
already presented in~\eqref{eq:darcy}, which we rewrite using the subscript $\Omega$ to identify quantities in $\Omega\setminus\gamma$
\begin{subequations}\label{eq:fracture_darcy}
\begin{align}\label{eq:fracture_darcy_domain}
    \begin{aligned}
        &\begin{aligned}
            & \bm{u}_\Omega + K_\Omega \nabla p_\Omega = \bm{0} \\
            &\nabla \cdot \bm{u}_\Omega = f_\Omega
        \end{aligned}
        && \text{in } \Omega\setminus\gamma,\\
        &\bm{u}_\Omega \cdot \bm{n}_\partial = \overline{u}_{\Omega} &&
        \text{on } \partial_u \Omega\setminus\partial\gamma ,\\
        & p_\Omega = \overline{p}_\Omega && \text{on } \partial_p {\Omega}\setminus\partial\gamma.
    \end{aligned}
\end{align}
We assume that also the flow in the fracture is governed by Darcy's law, however the differential operators operate now on the tangent
space. Yet, for the sake of simplicity, with an abuse of notation we use the same symbols to denote them.
The system of equations in the fracture is then given by
\begin{align}\label{eq:fracture_darcy_fracture}
    \begin{aligned}
        &\begin{aligned}
            &\epsilon_\gamma^{-1} \bm{u}_\gamma + K_\gamma \nabla p_\gamma = \bm{0} \\
            &\nabla \cdot \bm{u}_\gamma - \tr \bm{u}_{+} \cdot
            \bm{n} +
            \tr \bm{u}_{-} \cdot \bm{n} = f_\gamma
        \end{aligned}
        && \text{in } \gamma,\\
        &\bm{u}_\gamma \cdot \bm{n}_\partial = \overline{u}_\gamma &&
        \text{on } \partial_u \gamma,\\
        & p_\gamma = \overline{p}_\gamma && \text{on }
        \partial_p \gamma.
    \end{aligned}
\end{align}
Here, $\overline{u}_\gamma$ and $\overline{p}_\gamma$ are given boundary data, and we recall that possible fracture tips are
in $\partial_u \gamma$ with $\overline{u}_\gamma=0$. The parameter $K_\gamma$ is the tangential effective permeability in $\gamma$.
In the 2D setting, where the reduced fracture model is one-dimensional, $K_\gamma$ is a positive quantity. In the 3D setting, it may be in general
a rank-2 symmetric and positive tensor.
We may note in the equation representing the conservation of mass  the presence of an additional term that describes the flux exchange with the surrounding porous media. To close the problem we need to
complete the coupling between fracture and bulk, and we consider the following Robin-type condition on both
sides of $\gamma$
\begin{align}\label{eq:fracture_darcy_coupling}
    \epsilon_\gamma \tr \bm{u}_{\pm} \cdot \bm{n} \pm \kappa_\gamma (p_\gamma - \tr
    p_{\pm} )= 0 \quad \text{on } \gamma^{\pm},
\end{align}
with $\kappa_\gamma>0$ being the normal effective permeability. Problem~\eqref{eq:fracture_darcy} consists of the system of
equations that describe the Darcy velocity and pressure in both the fracture and
surrounding porous medium. An analysis may be found, for instance, in~\cite{Scotti2017} or \cite{Fumagalli2016a}.

\end{subequations}

The case of $N>1$ non-intersecting fractures the problem is analogous to the one just described where $\gamma=\cup_{i=1}^{N}\gamma_i$.
However, if two or more fractures intersect  we need to introduce new
conditions to describe the flux interchange between connected fractures. At each intersection $\iota$ we denote
with $I_\iota$ the set of intersecting fractures and we consider the following conditions on $\iota$,
\begin{gather}\label{eq:fracture_darcy_intersection}
    \begin{cases}
        \epsilon_\iota \alpha_j \tr \bm{u}_{j} \cdot \bm{t}_{j} + \kappa_\iota
        (p_\iota - \tr p_{j}) = 0\quad \forall \gamma_j\in I_\iota\\
        \displaystyle \sum_{\gamma_j\in I_\iota}\alpha_j \tr \bm{u}_{j} \cdot \bm{t}_{j}=0
    \end{cases}
    \quad \text{on } \iota,
\end{gather}
where $\epsilon_\iota$ is the measure of the
intersection, $p_\iota$ is the pressure at the
intersection, $\kappa_\iota$ is the permeability at the intersection
and $\alpha\in\{-1,1\}$ depends on the orientation chosen for the normal $\bm{t}_{j}$ to $\partial\gamma_j$ at the intersection.
Note that $\bm{t}_{j}$ is indeed on the tangent plane of $\gamma_j$. System~\eqref{eq:fracture_darcy_intersection} can be simplified by
noting that it implies that $p_\iota$ is equal to the average of the $p_{j}$.

\section{Weak formulation}\label{sec:weak_form}

The numerical scheme that we will present in Section \ref{sec:vem} is based on the weak
formulation of problem \eqref{eq:darcy}
and \eqref{eq:fracture_darcy}. Therefore, we will present in the following the functional setting and the weak form we have used as basis for the numerical discretization.
We indicate with $L^2(A)$ the Lebesgue space of square integrable functions on $A$, while
$H_{\operatorname{div}}(A)$ is the space of square integrable vector functions
whose distributional divergence is in $L^2(A)$. They are Hilbert spaces with standard norms and inner products.
In particular, we denote with $(\cdot, \cdot)_A$ the $L^2(A)$-scalar product. Moreover, given a functional
space $V$ and its dual $V^\prime$ we use $\langle a, b \rangle$, with $a\in V$ and $b\in V^\prime$ to denote the
duality pairing between the given functional spaces.

\subsection{Single-phase bulk flow without fractures}
If fractures are not present, the setting is rather standard.
For simplicity, we assume homogeneous essential boundary conditions
$\overline{u}_\Omega = 0$, otherwise a
lifting technique can be used to recover the original problem. We introduce
the following functional spaces for vector and scalar field, respectively,
\begin{gather}\label{eq:func_space_darcy}
    V(\Omega) = \left\{ \bm{v} \in H_{\operatorname{div}}(\Omega): \, \tr \bm{v} \cdot
    \bm{n}_\partial = 0 \text{ on } \partial_u\Omega \right\}
    \quad \text{and} \quad
    Q(\Omega) = L^2(\Omega).
\end{gather}
Here $\tr$ is the normal trace operator $\tr: H_{\operatorname{div}} (\Omega) \to
H^{-\frac{1}{2}}(\partial_u \Omega)$, which is linear and bounded, see~\cite{Boffi2013}.

We can now introduce the following bilinear forms and functionals
\begin{gather*}
    a_\Omega: V(\Omega) \times V(\Omega) \rightarrow \mathbb{R}: \quad
    a_\Omega(\bm{u}_\Omega,
    \bm{v}_\Omega) = (H_\Omega \bm{u}_\Omega, \bm{v}_\Omega)_\Omega\\
    b_\Omega: V(\Omega) \times Q(\Omega) \rightarrow \mathbb{R}: \quad b_\Omega(\bm{v}_\Omega, p_\Omega)
    = -(\nabla \cdot
    \bm{v}_\Omega, q_\Omega)_\Omega\\
    G_\Omega: V(\Omega) \rightarrow \mathbb{R}: \quad G_\Omega(\bm{v}_\Omega) = -\langle \tr \bm{v}_\Omega \cdot
    \bm{n}_\partial, \overline{p}_\Omega \rangle\\
    F_\Omega: Q(\Omega) \rightarrow \mathbb{R}: \quad F_\Omega(q_\Omega) = -(f_\Omega, q_\Omega)_\Omega
\end{gather*}
where $H_\Omega = K_\Omega^{-1}$. We assume that $K_\Omega \in [L^\infty(\Omega)]^{2 \times 2}$,
with $\underline{\alpha}||\bm{y}||^2\le \bm{y}^T K \bm{y}\le \overline{\alpha}||\bm{y}||^2$,
a.e. in $\Omega$, where $\bm{y}\in\mathbb{R}^n$ and $0<\underline{\alpha}\le \overline{\alpha}$.

Furthermore, we take $\overline{p}_\Omega
\in H^{\frac{1}{2}}_{00}(\partial_p \Omega)$, and $f_\Omega \in L^2(\Omega)$. Let us
note that $a_\Omega:\ V(\Omega)\times V(\Omega)\to \mathbb{R}$ is continuous, coercive and symmetric, being $K_\Omega$ symmetric.

We can now state the weak formulation of our problem: find $(\bm{u}_\Omega, p_\Omega) \in V(\Omega) \times
Q(\Omega)$ such that
\begin{gather}\label{eq:darcy_weak}
    \begin{aligned}
        & a_\Omega(\bm{u}_\Omega, \bm{v}_\Omega) + b_\Omega(\bm{v}_\Omega,
        p_\Omega) = G(\bm{v}_\Omega) && \forall \bm{v}_\Omega  \in V(\Omega)\\
        & b_\Omega(\bm{u}_\Omega, q_\Omega) = F_\Omega(q_\Omega) && \forall
        q_\Omega \in Q(\Omega)
    \end{aligned}.
\end{gather}
The previous problem is well posed, provided $\vert \partial_p
\Omega\vert >0$. See, for example, \cite{Boffi2013} for a proof.

\subsection{Fracture flow}

We extend now the weak formulation for problem \eqref{eq:fracture_darcy}, with
the simplifying assumption that only one fracture is considered. Its extension to multiple
fractures is straightforward, see for example \cite{Scotti2017,Boon2018}. Also in this case we assume homogeneous essential boundary
conditions, otherwise a lifting technique can be used.

We need to introduce the space $H_{\operatorname{div}}(\Omega\setminus \gamma)$ as the space of vector function
in $L^2(\Omega\setminus\gamma)$ (which may be identified by $L^2(\Omega)$ since $\gamma$ has zero measure) whose distributional divergence is in $L^2(\Sigma)$ for
all measurable $\Sigma\subset(\Omega\setminus\gamma)$. We need also to impose some extra regularity on the trace on $\gamma^{\pm}$, due to the Robin-type condition
\eqref{eq:fracture_darcy_coupling}. The reader may refer to
\cite{Martin2005,Frih2011,Boffi2013} for a
more detailed discussion on this matter.  In particular, we require that, for a $\bm{v}_{\Omega}\in H_{\operatorname{div}}(\Omega\setminus \gamma)$ ,
$\tr\bm{v}_{+}\cdot\bm{n} \in L^2(\gamma)$ and
$\tr\bm{v}_{-}\cdot\bm{n} \in L^2(\gamma)$,
where $\tr$ here indicates the trace of $\bm{v}$ on the two sides of the fracture. This space is equipped with the inner product
\begin{gather*}
    (\bm{u},\bm{v})_{H_{\operatorname{div}}(\Omega\setminus \gamma)}=
    (\bm{u},\bm{v})_{\Omega}+(\nabla\cdot\bm{u},\nabla\cdot\bm{v})_{\Omega}+
    (\tr\bm{u}_{+}\cdot\bm{n}, \tr\bm{v}_{+}\cdot\bm{n})_\gamma+(\tr\bm{u}_{-}\cdot\bm{n}, \tr\bm{v}_{-}\cdot\bm{n})_\gamma,
\end{gather*}
and induced norm.
The new space for vector fields in the bulk is
given by
\begin{gather*}
    \hat{V}(\Omega) = \left\{ \bm{v}_\Omega \in H_{\operatorname{div}}(\Omega\setminus\gamma): \, \tr
  \bm{v}_{\Omega} \cdot \bm{n}_\partial=0 \text{ on } \partial_u\Omega \right\}.
\end{gather*}
The functional spaces for vector and
scalar fields defined on the fracture are
\begin{gather*}
    V(\gamma) = \left\{ \bm{v}_\gamma \in H_{\operatorname{div}}(\gamma): \, \tr
    \bm{v}_\gamma \cdot
    \bm{n}_\partial = 0
    \right\}
    \quad \text{and} \quad
    Q(\gamma) = L^2(\gamma),
\end{gather*}
where in this case the trace operator in $V(\gamma)$ is given by $\tr:
H_{\operatorname{div{}}}(\gamma) \rightarrow H^{-\frac{1}{2}}(\partial_u
\gamma)$. Note that in the case of 2D problems like the ones treated in this
work, $V(\gamma)$ is in fact a subspace of $H^1(\gamma)$ and
the trace reduces to the value at the boundary point.

We introduce now the bilinear forms and functional for the weak formulation of
problem \eqref{eq:fracture_darcy}.
First, we modify the bilinear form $a_\Omega$
by taking into account the coupling terms from
\eqref{eq:fracture_darcy_coupling} as
\begin{gather*}
    \hat{a}_\Omega: \hat{V}(\Omega) \times \hat{V}(\Omega)
    \rightarrow \mathbb{R}: \quad
    \hat{a}_\Omega(\bm{u}_\Omega,
    \bm{v}_\Omega) = a_\Omega(\bm{u}_\Omega, \bm{v}_\Omega)_\Omega +
    \sum_{*\in\{+, -\}}(\eta_\gamma \tr \bm{u}_* \cdot \bm{n},
    \tr \bm{u}_{*} \cdot \bm{n} )_{\gamma}
\end{gather*}
where $\eta_\gamma = \epsilon_\gamma \kappa_\gamma^{-1}$ and we have assumed
that $\eta_\gamma \in L^\infty(\gamma)$.
Second, the bilinear forms associated
with the fracture are given by
\begin{gather*}
    a_\gamma: V(\gamma) \times V(\gamma) \rightarrow \mathbb{R}: \quad
    a_\gamma(\bm{u}_\gamma, \bm{v}_\gamma) = (H_\gamma \bm{u}_\gamma,
    \bm{v}_\gamma )_\gamma\\
    b_\gamma: V(\gamma) \times Q(\gamma) \rightarrow \mathbb{R}: \quad
    b_\gamma(\bm{v}_\gamma, p_\gamma) = -(\nabla \cdot \bm{v}_\gamma,
    p_\gamma)_\gamma\\
    G_\gamma: V(\gamma) \rightarrow \mathbb{R}: \quad
    G_\gamma(\bm{v}_\gamma) = - \langle \tr \bm{v}_\gamma \cdot \bm{n}_\partial,
    \overline{p}_\gamma \rangle\\
    F_\gamma: V(\gamma) \times \mathbb{R}: \quad
    F_\gamma(q_\gamma) = -(f_\gamma, q_\gamma)_\gamma
\end{gather*}
where we have $H_\gamma^{-1} = \epsilon_\gamma K_\gamma$ and we have assumed
that $H_\gamma \in L^\infty(\gamma)$, $\overline{p}_\gamma \in
H^{\frac{1}{2}}(\partial_p \gamma)$, and $f_\gamma \in L^2(\gamma)$. Third, we
introduce the bilinear forms responsible for the flux exchange between the
fracture and the bulk medium
\begin{gather*}
    c^\pm: \hat{V}(\Omega) \times Q(\gamma) \rightarrow \mathbb{R}: \quad
    c^\pm(\bm{u}_\Omega, q_\gamma) = \pm (\tr \bm{u}_{\pm} \cdot \bm{n},
    q_\gamma )_\gamma\\
    c: \hat{V}(\Omega) \times Q(\gamma) \rightarrow \mathbb{R}: \quad
    c(\bm{u}_\Omega, q_\gamma) = \sum_{* \in \{+, -\}} c^*(\bm{u}_\Omega, q_\gamma).
\end{gather*}
We finally can write the weak formulation for problem
\eqref{eq:fracture_darcy}: find $(\bm{u}_\Omega, p_\Omega,
\bm{u}_\gamma, p_\gamma) \in \hat{V}(\Omega) \times Q(\Omega) \times V(\gamma)
\times Q(\gamma)$ such that
\begin{gather}\label{eq:weak_fracture}
    \begin{aligned}
        & \hat{a}_\Omega(\bm{u}_\Omega, \bm{v}_\Omega) + b_\Omega(\bm{v}_\Omega,
        p_\Omega) + c(\bm{v}_{\Omega}, p_\gamma) =
        G_\Omega(\bm{v}_\Omega) && \forall \bm{v}_\Omega \in \hat{V}(\Omega)\\
        & b_\Omega(\bm{u}_\Omega, q_\Omega) = F_\Omega(q_\Omega) && \forall
        q_\Omega \in Q(\Omega)\\
        &a_\gamma(\bm{u}_\gamma, \bm{v}_\gamma) + b_\gamma(\bm{v}_\gamma,
        p_\gamma) = G_\gamma(\bm{v}_\gamma) && \forall \bm{v}_\gamma \in
        V(\gamma)\\
        &b_\gamma(\bm{u}_\gamma, q_\gamma) + c(\bm{u}_\Omega, q_\gamma) =
        F_\gamma(q_\gamma) && \forall q_\gamma \in Q(\gamma)
    \end{aligned}
\end{gather}
The reader can refer to \cite{DAngelo2011,Formaggia2012,DelPra2015a} for proofs of the
well posedness of the problem, provided suitable boundary conditions.

\section{Numerical approximation by MVEM}\label{sec:vem}

The challenges in terms of heterogeneity of physical data and complexity of the geometry
due to the presence of fractures influence the choice of the numerical scheme.
A possible choice is the mixed finite element method, see
\cite{Raviart1977,Roberts1991,Boffi2013}. However,
this class of methods, capable of providing accurate results for pressure and
velocity fields, even in the presence of high heterogeneities, requires grids made either of simplexes (triangles of tetrahedra)
or quad/hexahedra. This can be inefficient for the problem at hand, where instead methods able to operate on grids formed
by  arbitrary polytopes are rather appealing. For this reason finite volume schemes, see~\cite{Droniou2013} for a review, are very much used
in practice. However, they normally treat the primal formulation and require good quality grids to obtain an accurate solution and a good reconstruction of the velocity field. Indeed, it is known that convergence of the method is guaranteed only if the grid has special properties.

Therefore, we focus here our attention on
the low-order Mixed Virtual Element Method, a numerical schemes that operates on polytopal grids and that has shown to be rather robust with respect to irregularities in the data and in the computational grid. We consider first the case of  porous medium without fractures, focusing on
problems with highly heterogeneous permeability, and then the case of a fractured porous medium, using the model described in Subsection
\ref{subsec:vem_frac}.

The actual implementation in PorePy adopts a flux mortar technique that allows
non-conforming coupling between inter-dimensional grids. We do not exploit the
possibility of having grids non-conforming to the fractures in this work, nevertheless in Subsection \ref{subsec:vem_frac} we will describe the mortar approach more in detail.

\subsection{Bulk flow without fractures}\label{subsec:vem_single}

In this part we present the MVEM discretization of problem
\eqref{eq:darcy_weak}. A key point of the virtual method is to use an implicit definition of suitable
basis functions, and obtain
computable discrete local matrices by manipulating the different terms in the weak formulation appropriately.
In this work we consider only the low order case, yet the method can be extended to higher order formulations.

We indicate the computational grid, approximation of $\Omega$, as
$\mathcal{T}(\Omega)$. We assume that $\Omega$ has polygonal boundary, so that
$\mathcal{T}(\Omega)$ covers $\Omega$ exactly.
The set of faces of $\mathcal{T}(\Omega)$ is denoted as $\mathcal{E}(\Omega)$, with the
distinction between the internal and boundary faces indicated by
$\mathcal{E}({\interior{\Omega}})$ and
$\mathcal{E}({\partial \Omega}$), respectively. We also specify the edges on a
specific portion of the boundary of $\Omega$ as $\mathcal{E}(\partial_u \Omega)$
and $\mathcal{E}(\partial_p \Omega)$.
We clearly have $\mathcal{E}(\interior{\Omega}) \cup
\mathcal{E}(\partial \Omega) = \mathcal{E}(\Omega)$ as well as
$\mathcal{E}(\interior{\Omega}) \cap
\mathcal{E}(\partial \Omega) = \emptyset$.
In the sequel, we generally indicate as $C \in \mathcal{T}(\Omega)$ a grid cell and $e \in
\mathcal{E}(\Omega)$ a face between two cells. Element $C$ can be a generic polygon (polyhedra in the 3D case).

We introduce the finite dimensional subspaces, approximation of the continuous
spaces given in \eqref{eq:func_space_darcy}, as
\begin{gather*}
    V_{h}(\Omega) = \left\{ \bm{v}_\Omega \in V(\Omega): \, \nabla \cdot
    \bm{v}_\Omega|_C \in \mathbb{P}_0(C) \text{ and } \nabla \times
    \bm{v}_\Omega|_C = \bm{0},\,
    \forall C \in \mathcal{T}(\Omega), \right. \\
    \left. \, \tr \bm{v}_\Omega \cdot \bm{n}_e \in
    \mathbb{P}_0(e), \, \forall e \in \mathcal{E}(\Omega) \right\},
\end{gather*}
with $\mathbb{P}_0(X)$ being the space of constant polynomials on $X$, while $\tr$ and $\bm{n}_e$
the trace and the normal associated to edge $e$. For the scalar
field we set
\begin{gather*}
    Q_h({\Omega}) = \left\{ q_\Omega \in Q(\Omega): \, q_\Omega |_C \in
    \mathbb{P}_0(C) \, \forall C \in \mathcal{T}(\Omega)\right\}.
\end{gather*}

Clearly $V_h({\Omega}) \subset V(\Omega)$ and $Q_h({\Omega}) \subset
Q(\Omega)$. The degrees of freedom associated with $V_h({\Omega})$ and $Q_h({\Omega})$ are one scalar
value for each face and one scalar value for each cell, respectively.
More precisely, if we generically indicate with $\dof_i$ the functional associated with the
$i$-th degree of freedom, we have, for a $\bm{v}_\Omega\in V_h(\Omega)$ and a $q_\Omega \in Q_h({\Omega})$
\begin{gather*}
    \dof_i \bm{v}_\Omega = \tr\bm{v}_\Omega\cdot \bm{n}_{e_i}
    \quad\text{and}\quad
    \dof_i  q_\Omega = q_\Omega\vert_{C_i},
\end{gather*}
where $e_i$ and $C_i$ are the $i$-th edge and cell, respectively, and $\tr$ now
indicates the trace associated to the edge $e_i$.

Moreover, we
can observe that in case of triangular grids $V_h({\Omega})$ coincides with $\mathbb{RT}_0(\Omega)$, so the former can be viewed
as a generalization of the well known Raviart-Thomas finite elements.

By performing exact integration, the numerical approximation of the
bilinear form $b_\Omega$ and of the functionals
$G_\Omega$, $F_\Omega$ are computable with the
given definition of the discrete spaces. However, for the term ${a}_\Omega$
we need further manipulations to obtain a computable expression.
To this purpose, we define a suitable subspace of $V_h(\Omega)$, defined as
\begin{gather*}
    \mathcal{V}_h(\Omega) = \left\{ \bm{v}_\Omega \in V_h(\Omega): \,
    \bm{v}_\Omega |_C = K_C \nabla v_C \text{ for a } v_C \in
    \mathbb{P}_1(C) \forall C \in \mathcal{T}(\Omega)\right\},
\end{gather*}
where $K_C$ is a suitable constant approximation of $K_\Omega\vert_C$, and
we define a projection operator $\Pi_\Omega: V_h(\Omega) \rightarrow
\mathcal{V}_h(\Omega)$ so that for a $\bm{v} \in V_h(\Omega)$ we
have
\begin{gather*}
a_\Omega(\bm{v} - \Pi_\Omega \bm{v}, \bm{w}) =0,\quad \forall \bm{w} \in \mathcal{V}_h(\Omega).
\end{gather*}
We now set
$T_\Omega = I - \Pi_\Omega$, where $T_\Omega: V_h(\Omega) \rightarrow
\mathcal{V}_h^\perp(\Omega)$ and the orthogonality condition is governed by the
 bilinear form $a_\Omega$, which, being symmetric, continuous and coercive, defines an
 inner product. Indeed, from the definition of $\Pi_\Omega$ we have
$a_\Omega(T_\Omega \bm{v}_\Omega, \Pi_\Omega \bm{w}_\Omega) = 0$ for all
$\bm{v}_\Omega, \bm{w}_\Omega \in V_h(\Omega)$. Considering this fact,
we have the following decomposition
\begin{gather*}
    a_\Omega ( \bm{u}_\Omega, \bm{v}_\Omega ) = a_\Omega ( (\Pi_\Omega +
    T_\Omega) \bm{u}_\Omega, (\Pi_\Omega + T_\Omega) \bm{v}_\Omega ) =
    a_\Omega ( \Pi_\Omega
    \bm{u}_\Omega, \Pi_\Omega \bm{v}_\Omega ) +
    a_\Omega ( T_\Omega \bm{u}_\Omega, T_\Omega \bm{v}_\Omega ).
\end{gather*}
Now, thanks to the definition of $\mathcal{V}_h(\Omega)$ the first term is
computable in terms of the degrees of freedom, see for instance \cite{Fumagalli2016b}, but not the second one. However, since it gives the contribution of
$a_\Omega$ only on $\mathcal{V}_h^\perp(\Omega)$, it can be approximated with a suitable stabilizing bilinear form
$s:V_h(\Omega) \times V_h(\Omega) \rightarrow \mathbb{R}$, i.e.
\begin{gather*}
    a_\Omega ( T_\Omega \bm{u}_\Omega, T_\Omega \bm{v}_\Omega ) \approx
    s_\Omega( \bm{u}_\Omega, \bm{v}_\Omega).
\end{gather*}
For more details about $s_\Omega$ refer to the works
\cite{Brezzi2014,BeiraodaVeiga2014b,BeiraoVeiga2016,Fumagalli2016a,Fumagalli2017a,Dassi2019}.
The form $s_\Omega$ must satisfy the following equivalence condition:
\begin{gather*}
    \exists \upsilon_*, \upsilon^* \in \mathbb{R}^+: \quad
    \upsilon_* a_\Omega( \bm{u}_\Omega, \bm{v}_\Omega) \leq
    s_\Omega( \bm{u}_\Omega, \bm{v}_\Omega) \leq \upsilon^*
    a_\Omega( \bm{u}_\Omega, \bm{v}_\Omega) \quad \forall \bm{u}_\Omega, \bm{v}_\Omega
    \in V_h(\Omega).
\end{gather*}
The illustrate our choice of $s_\Omega$, let us denote with
$\bm{\varphi}$ a generic element of the basis of $V_h(\Omega)$.
The stabilization term, in our case, can be computed as
\begin{gather}\label{eq:stabilization_mvem}
    s_\Omega (\bm{\varphi}_{ \theta}, \bm{\varphi}_{\chi}) =\sum_{C \in
    \mathcal{T}(\Omega)}
    \norm{H_\Omega}_{L^\infty(C)}
    \sum_{i=1}^{N_{dof}(C)} \dof_i(T_\Omega \bm{\varphi}_{\theta}) \dof_i
    (T_\Omega \bm{\varphi}_{\chi}),
\end{gather}
where $N_{dof}(C)$ is the total number of degrees of freedom for the vector
field for the cell $C$
and $\dof_i$ gives the value of the argument at the $i^{th}$-dof. The $K_\Omega$
norm is a scaling factor in order to consider also strong oscillations of
physical parameters. With the definition of the stabilization term now all the terms are computable and the global system can be assembled.
For more details on the actual computation of the local matrices refer to
\cite{Fumagalli2016a,BeiraodaVeiga2014a}.

\subsection{Fracture flow}\label{subsec:vem_frac}

We introduce now the numerical scheme used for the approximation of problem
\eqref{eq:weak_fracture}. We consider the notations and terms for the porous
media from the previous section. In fact the derivation of the discrete setting
for the porous media is similar to what already presented. We focus now on the
fracture discretization as well as on the coupling term with the surrounding porous
media.

In particular, for the implementation we have chosen PorePy
\cite{Keilegavlen2019a}, that  considers an additional interface $\gamma^\pm$ between the
fracture and the porous media along with a flux mortar technique to couple
domains of different dimensions, allowing also non-conforming grids between the
domains. However, to avoid additional complexity we consider only conforming
grids so that the mortar variable behaves as a Lagrange multiplier $\lambda_h$. The latter is
the normal flux exchange from the higher to lower dimensional domain. See Figure
\ref{fig:vem_frac} as an example. Geometrically \textit{i)} the interface between the porous
media and the fracture, \textit{ii)} the fracture, and  \textit{iii)} the two interfaces coincide but they
are represented by different objects with suitable operators for their coupling.
In the case of conforming discretizations these operators simply map the
corresponding degrees of freedom, however in the case of non-conforming
discretizations projection operators should be considered.
\begin{figure}[tbp]
    \centering
    \subfloat[Interfaces\label{fig:vem_frac}]
    {
\begingroup%
  \makeatletter%
  \providecommand\color[2][]{%
    \errmessage{(Inkscape) Color is used for the text in Inkscape, but the package 'color.sty' is not loaded}%
    \renewcommand\color[2][]{}%
  }%
  \providecommand\transparent[1]{%
    \errmessage{(Inkscape) Transparency is used (non-zero) for the text in Inkscape, but the package 'transparent.sty' is not loaded}%
    \renewcommand\transparent[1]{}%
  }%
  \providecommand\rotatebox[2]{#2}%
  \newcommand*\fsize{\dimexpr\f@size pt\relax}%
  \newcommand*\lineheight[1]{\fontsize{\fsize}{#1\fsize}\selectfont}%
  \ifx\svgwidth\undefined%
    \setlength{\unitlength}{154.02376791bp}%
    \ifx\svgscale\undefined%
      \relax%
    \else%
      \setlength{\unitlength}{\unitlength * \real{\svgscale}}%
    \fi%
  \else%
    \setlength{\unitlength}{\svgwidth}%
  \fi%
  \global\let\svgwidth\undefined%
  \global\let\svgscale\undefined%
  \makeatother%
  \begin{picture}(1,0.56888735)%
    \lineheight{1}%
    \setlength\tabcolsep{0pt}%
    \put(0,0){\includegraphics[width=\unitlength,page=1]{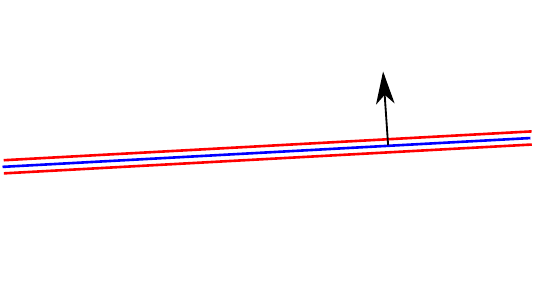}}%
    \put(0.73698117,0.3700154){\color[rgb]{0,0,0}\makebox(0,0)[lt]{\lineheight{1.25}\smash{\begin{tabular}[t]{l}$\bm{n}$\end{tabular}}}}%
    \put(0.38240354,0.04424196){\color[rgb]{0,0,0}\makebox(0,0)[lt]{\lineheight{1.25}\smash{\begin{tabular}[t]{l}$\Omega^-$\end{tabular}}}}%
    \put(0,0){\includegraphics[width=\unitlength,page=2]{vem_frac.pdf}}%
    \put(0.39773586,0.49800283){\color[rgb]{0,0,0}\makebox(0,0)[lt]{\lineheight{1.25}\smash{\begin{tabular}[t]{l}$\Omega^+$\end{tabular}}}}%
    \put(0.01143264,0.27768759){\color[rgb]{0,0,0}\makebox(0,0)[lt]{\lineheight{1.25}\smash{\begin{tabular}[t]{l}$\gamma$\end{tabular}}}}%
    \put(0.45446302,0.21271128){\color[rgb]{0,0,0}\makebox(0,0)[lt]{\lineheight{1.25}\smash{\begin{tabular}[t]{l}$\gamma^-$\end{tabular}}}}%
    \put(0.45836869,0.31231744){\color[rgb]{0,0,0}\makebox(0,0)[lt]{\lineheight{1.25}\smash{\begin{tabular}[t]{l}$\gamma^+$\end{tabular}}}}%
  \end{picture}%
\endgroup%
}
    \hspace*{0.1\textwidth}
    \subfloat[Degrees of freedom \label{fig:vem_frac_dof}]
    {
\begingroup%
  \makeatletter%
  \providecommand\color[2][]{%
    \errmessage{(Inkscape) Color is used for the text in Inkscape, but the package 'color.sty' is not loaded}%
    \renewcommand\color[2][]{}%
  }%
  \providecommand\transparent[1]{%
    \errmessage{(Inkscape) Transparency is used (non-zero) for the text in Inkscape, but the package 'transparent.sty' is not loaded}%
    \renewcommand\transparent[1]{}%
  }%
  \providecommand\rotatebox[2]{#2}%
  \newcommand*\fsize{\dimexpr\f@size pt\relax}%
  \newcommand*\lineheight[1]{\fontsize{\fsize}{#1\fsize}\selectfont}%
  \ifx\svgwidth\undefined%
    \setlength{\unitlength}{154.02376791bp}%
    \ifx\svgscale\undefined%
      \relax%
    \else%
      \setlength{\unitlength}{\unitlength * \real{\svgscale}}%
    \fi%
  \else%
    \setlength{\unitlength}{\svgwidth}%
  \fi%
  \global\let\svgwidth\undefined%
  \global\let\svgscale\undefined%
  \makeatother%
  \begin{picture}(1,0.56888735)%
    \lineheight{1}%
    \setlength\tabcolsep{0pt}%
    \put(0.64535401,0.04424196){\color[rgb]{0,0,0}\makebox(0,0)[lt]{\lineheight{1.25}\smash{\begin{tabular}[t]{l}$\Omega^-$\end{tabular}}}}%
    \put(0,0){\includegraphics[width=\unitlength,page=1]{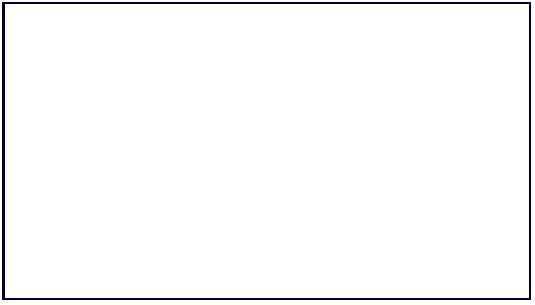}}%
    \put(0.66068633,0.49800283){\color[rgb]{0,0,0}\makebox(0,0)[lt]{\lineheight{1.25}\smash{\begin{tabular}[t]{l}$\Omega^+$\end{tabular}}}}%
    \put(0.50315756,0.22244991){\color[rgb]{0,0,0}\makebox(0,0)[lt]{\lineheight{1.25}\smash{\begin{tabular}[t]{l}$\gamma^-$\end{tabular}}}}%
    \put(0.50706323,0.32205606){\color[rgb]{0,0,0}\makebox(0,0)[lt]{\lineheight{1.25}\smash{\begin{tabular}[t]{l}$\gamma^+$\end{tabular}}}}%
    \put(0,0){\includegraphics[width=\unitlength,page=2]{vem_frac_dof.pdf}}%
  \end{picture}%
\endgroup%
}
    \caption{On the left, the hybrid-dimensional representation of a fracture immersed in a porous
    media with the two interfaces $\gamma^+$ and $\gamma^-$, in red. On the
    right, the representation of the degrees of freedom for vector fields.}
    \label{fig:vem_frac}
\end{figure}

As done before, we consider the special case of a single fracture, being
its generalization straightforward. First of all, the velocity degrees of freedom for the rock
matrix in the proximity of the fracture are doubled as Figure \ref{fig:vem_frac_dof} shows.
We
can thus represent $\tr \bm{u}_{\pm} \cdot \bm{n} = \lambda^\pm_h$ for both sides $\pm$
of the fracture itself. The term $\hat{a}_\Omega$
involves the actual integration of the basis functions for each grid cells,
which is not possible since they are not, in general, analytically known.

Many of the following steps are similar to what already done for the bulk porous
media. We introduce a tessellation of $\gamma$ into non-overlapping cells (segments in
this case), the grid is indicated with $\mathcal{T}(\gamma)$ and the set of
faces (edges) as $\mathcal{E}(\gamma)$. Also in this case, we divide the internal
faces and the external faces $\mathcal{E}(\interior{\gamma})$ and
$\mathcal{E}(\partial \gamma)$. Moreover, the latter can also be divided into
subset depending on the boundary conditions $\mathcal{E}(\partial_u \gamma)$ and
$\mathcal{E}(\partial_p \gamma)$. Clearly, we have $\mathcal{E}(\gamma) =
\mathcal{E}(\interior{\gamma}) \cup \mathcal{E}(\partial \gamma)$ as well as
$\mathcal{E}(\partial_u \gamma) \cup \mathcal{E}(\partial_p
\gamma)=\mathcal{E}(\partial \gamma)$.
We introduce the functional spaces for the variables
defined on the fracture, for the vector fields we have
\begin{gather*}
    V_{h}(\gamma) = \left\{ \bm{v}_\gamma \in V(\gamma): \, \nabla \cdot
    \bm{v}_\gamma|_C \in \mathbb{R} \text{ and } \nabla \times
    \bm{v}_\gamma|_C = \bm{0}\,
    \forall C \in \mathcal{T}(\gamma),  \, \tr \bm{v}_\gamma \cdot \bm{n}_e \in
    \mathbb{R} \, \forall e \in \mathcal{E}(\gamma) \right\},
\end{gather*}
while for the scalar fields we consider the discrete space
\begin{gather*}
    Q_h(\gamma) = \left\{ q_\gamma\in Q(\gamma): \, q_\gamma |_C \in
    \mathbb{R} \, \forall C \in \mathcal{T}(\gamma) \right\}.
\end{gather*}
By keeping the same approach as before, we assume exact integration so that the numerical approximation of the
bilinear form $b_\gamma$ as well as functionals
$G_\gamma$ and $F_\gamma$ are computable with the
given definition of the discrete spaces. The term $a_\gamma$ is not directly
computable, we thus introduce the subspace of $V_h(\gamma)$ as
\begin{gather*}
    \mathcal{V}_h(\gamma) = \left\{ \bm{v}_\gamma \in V_h(\gamma):\,
    \bm{v}_\gamma |_C= K_\gamma|_C \nabla v_C \text{ for a } v_C \in
    \mathbb{P}_1(C) \forall C \in \mathcal{T}(\gamma) \right\}.
\end{gather*}
We introduce the projection operator $\Pi_\gamma$ from $V_h(\gamma) \rightarrow
\mathcal{V}_h(\gamma)$ such that for a $\bm{v} \in V_h(\gamma)$ we have $a_\gamma(\bm{v} - \Pi_\gamma \bm{v}, \bm{w}) =
0$ for all $\bm{w} \in \mathcal{V}_h$. By introducing the operator $T_\gamma = I
- \Pi_\gamma$, we have the decomposition
\begin{gather*}
    a_\gamma(\bm{u}_\gamma, \bm{v}_\gamma)  =
    a_\gamma((\Pi_\gamma + T_\gamma)\bm{u}_\gamma, (\Pi_\gamma + T_\gamma)
    \bm{v}_\gamma)  =
    a_\gamma(\Pi_\gamma \bm{u}_\gamma, \Pi_\gamma
    \bm{v}_\gamma)+     a_\gamma(T_\gamma \bm{u}_\gamma, T_\gamma
    \bm{v}_\gamma)
\end{gather*}
By the definition of $\mathcal{V}_h(\gamma)$ the first term is now computable,
while the second term, which is not computable, is replaced by the stabilization
term
\begin{gather*}
    a_\gamma(T_\gamma \bm{u}_\gamma, T_\gamma \bm{v}_\gamma) \approx s_\gamma
    (\bm{u}_\gamma, \bm{v}_\gamma)
\end{gather*}
with the request that $s_\gamma$ scales as $a_\gamma$, meaning that
\begin{gather*}
    \exists \upsilon_*, \upsilon^* \in \mathbb{R}: \quad \upsilon_*
    s_\gamma(\bm{u}_\gamma, \bm{v}_\gamma) \leq a_\gamma(\bm{u}_\gamma,
    \bm{v}_\gamma) \leq \upsilon^* s_\gamma (\bm{u}_\gamma, \bm{v}_\gamma) \quad
    \forall \bm{u}_\gamma, \bm{v}_\gamma \in V_h(\gamma).
\end{gather*}
Denoting an element of the basis of $V_h(\gamma)$ as $\bm{\phi}$, the actual
construction of $s_\gamma$ is given by the formula
\begin{gather*}
    s_\gamma (\bm{\phi}_{ \theta}, \bm{\phi}_{\chi}) = \sum_{C \in
    \mathcal{T}(\gamma)}
    h\norm{K_\gamma^{-1}}_{L^\infty(C)}
    \sum_{i=1}^{N_{dof}(C)} \dof_i(T_\gamma \bm{\phi}_{\theta}) \dof_i
    (T_\gamma \bm{\phi}_{\chi}),
\end{gather*}
with $h$ the diameter of the current cell $C$. With the previous choices all the terms are computable and the fracture problem can be assembled.
For more details see \cite{Brezzi2014,BeiraodaVeiga2014b,BeiraoVeiga2016,Fumagalli2016a,Fumagalli2017a}.

\begin{figure}[tbp]
    \centering
\begingroup%
  \makeatletter%
  \providecommand\color[2][]{%
    \errmessage{(Inkscape) Color is used for the text in Inkscape, but the package 'color.sty' is not loaded}%
    \renewcommand\color[2][]{}%
  }%
  \providecommand\transparent[1]{%
    \errmessage{(Inkscape) Transparency is used (non-zero) for the text in Inkscape, but the package 'transparent.sty' is not loaded}%
    \renewcommand\transparent[1]{}%
  }%
  \providecommand\rotatebox[2]{#2}%
  \newcommand*\fsize{\dimexpr\f@size pt\relax}%
  \newcommand*\lineheight[1]{\fontsize{\fsize}{#1\fsize}\selectfont}%
  \ifx\svgwidth\undefined%
    \setlength{\unitlength}{154.02376791bp}%
    \ifx\svgscale\undefined%
      \relax%
    \else%
      \setlength{\unitlength}{\unitlength * \real{\svgscale}}%
    \fi%
  \else%
    \setlength{\unitlength}{\svgwidth}%
  \fi%
  \global\let\svgwidth\undefined%
  \global\let\svgscale\undefined%
  \makeatother%
  \begin{picture}(1,0.56888735)%
    \lineheight{1}%
    \setlength\tabcolsep{0pt}%
    \put(0,0){\includegraphics[width=\unitlength,page=1]{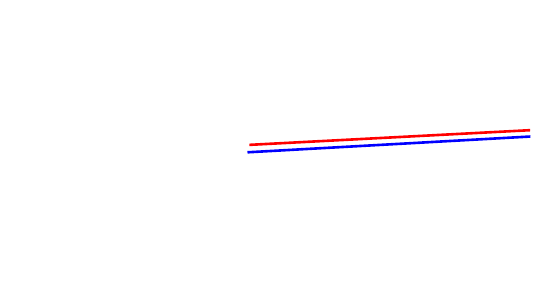}}%
    \put(0.64535401,0.04424196){\color[rgb]{0,0,0}\makebox(0,0)[lt]{\lineheight{1.25}\smash{\begin{tabular}[t]{l}$\Omega^-$\end{tabular}}}}%
    \put(0,0){\includegraphics[width=\unitlength,page=2]{vem_frac_lambda.pdf}}%
    \put(0.66068633,0.49800283){\color[rgb]{0,0,0}\makebox(0,0)[lt]{\lineheight{1.25}\smash{\begin{tabular}[t]{l}$\Omega^+$\end{tabular}}}}%
    \put(0.50315756,0.22244991){\color[rgb]{0,0,0}\makebox(0,0)[lt]{\lineheight{1.25}\smash{\begin{tabular}[t]{l}$\gamma^-$\end{tabular}}}}%
    \put(0.50706323,0.32205606){\color[rgb]{0,0,0}\makebox(0,0)[lt]{\lineheight{1.25}\smash{\begin{tabular}[t]{l}$\gamma^+$\end{tabular}}}}%
    \put(0,0){\includegraphics[width=\unitlength,page=3]{vem_frac_lambda.pdf}}%
    \put(0.86702197,0.30637187){\color[rgb]{0,0,0}\makebox(0,0)[lt]{\lineheight{1.25}\smash{\begin{tabular}[t]{l}$\gamma$\end{tabular}}}}%
  \end{picture}%
\endgroup%

    \caption{Representation of the conforming computational grids for the porous
    media, the fracture, and the two interfaces.}
    \label{fig:vem_frac_lambda}
\end{figure}

To couple the bulk and fracture flow, a Lagrange multiplier $\lambda^\pm_h$ is
used to represent the flux exchange between the fracture and the surrounding
porous media. We assume conforming grids, meaning that the
fracture grid is conforming with the interface grid as well as the faces of the
porous media are conforming with the interface grid. See Figure \ref{fig:vem_frac_lambda} as an
example. For space compatibility, we assume the Lagrange multiplier be a
piece-wise constant polynomial.
The interface condition \eqref{eq:fracture_darcy_coupling} is directly computable
with the degrees of freedom introduced providing a suitable projection of the
pressure $p_\Omega$ at the fracture interface. Our choice is to consider the
same value of the pressure at neighbouring cells, however other approaches can be
used, see for example \cite{Roberts1991}.

\section{Grid generation}\label{sec:grid_generation}

The generation of grids for realistic fractured porous media geometries is a
challenging task, whose complete automatic solution is still an open problem,
particularly for 3D configurations. We here give a brief overview of some
techniques that have been proposed, with no pretense of being exhaustive.

\subsection{Constrained Delaunay}

The generation of a grid of simplexes (triangles in 2D, tetrahedra in 3D)
conformal to a fracture network may be performed in principle by employing a
constrained Delaunay algorithm. It is an extension of the well known Delaunay
algorithm to the case where the mesh has to honor internal constraints (or
describe a non-convex domain). Usually it starts from a representation of
the domain and in 3D it first generates constrained Delaunay triangulation on
the fracture and boundary geometry, then new nodes are added in the domain to
generate a final grid that satisfies a relaxed Delaunay criterion to honour the
internal interfaces. The description of the constrained Delaunay procedure may
be found, for instance, in~\cite{Cheng2012}. Another general reference for mesh
generation procedures is~\cite{Frey2013}. However, in practical situations
several issues may arise. The presence of fractures intersecting with small
angles, for instance, may produce an excessive refinement near the intersections
in order to maintain the Delaunay property. In 3D there is the additional issue
of the possible generation of extremely badly distorted elements, often called
slivers, whose automatic removal is problematic, when not impossible, under the
constraint of  conformity with complex internal surfaces.

Several techniques have been proposed to ameliorate the procedure. For instance
in \cite{Mustapha2007} and \cite{Mustapha2014} the authors present a procedure
that modifies the fracture network trying to maintain its characteristics of
connectivity and effective permeability, while eliminating geometrical
situations where that may impair the effectiveness of a Delaunay triangulation.
In the second reference, a special decision strategy (called ``Gabriel
criterion'') is used to select a part of the fracture network to which
triangulation can be constrained without leading to an excessive degradation in
cells quality, or excessively fine grids. The procedure has proved rather
effective on moderately complex network in 2D, while the results for 3D
configurations seem less convincing.

We mention for completeness that an alternative procedure for generating
simplicial grids is the one based on the front advancing technique (maybe
coupled with the Delaunay procedure). It is implemented in several software
tools, see for instance~\cite{Schoeberl1997,Frey1998}. However, its use in the
context of fractures media is at the moment very limited, probably because of
the lack of results of the termination of the procedure, contrary to the
Delaunay algorithm where one can prove that, under mild conditions, the
generation terminates in a finite number of steps. Moreover it has a much higher
computational cost. The interested reader may consult the cited references.

In our case PorePy considers the software Gmsh \cite{Geuzaine2009} for
the generation of the Delaunay bidimensional grids. The grid size in the
configuration file is specifically tuned to obtain high quality triangles. Indeed, we consider distances between fractures, between a fracture and the domain
boundary, and length of fracture branches. With these precautions, we usually
obtain quality grids that are suitable for numerical studies.

\subsection{Grids cut by the fracture network}\label{cart} 
\begin{figure}[htb]
    \centering
    \includegraphics[width=0.8\textwidth]{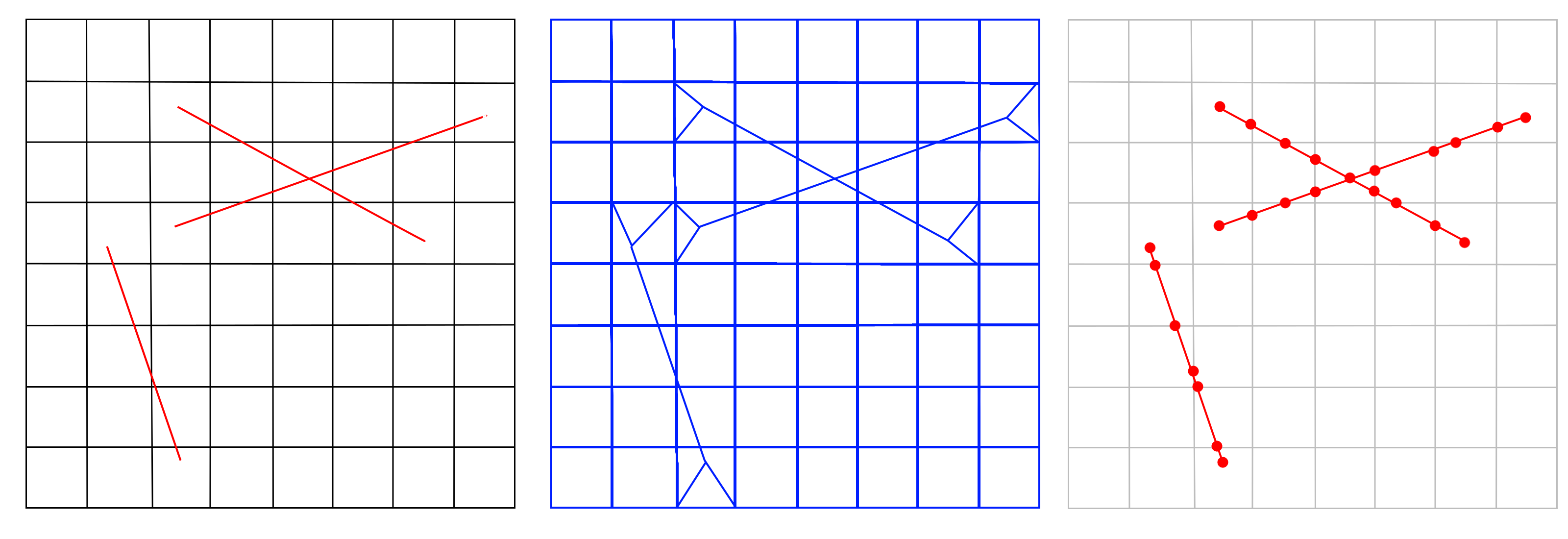}
    \caption{Creation of a polygonal mesh from a regular Cartesian grid.}
    \label{cut}
\end{figure}

An alternative
possibility to create a grid conforming to fractures or, in general, planar
interfaces, consists in cutting a regular Cartesian or simplex mesh, as shown in
Figure \ref{cut} for the case of a Cartesian mesh. The resulting grid will be
formed by polytopal elements in the vicinity of the fractures. The main issue
in this procedure is the possible generation of badly shaped or very small
elements as a consequence of the cut. Another technical problem is the necessity
of having efficient techniques for computing intersections and constructing the
polytope. To this respect, one may adopt the tools available in specialized
libraries like CGAL~\cite{cgalweb}, or developed by the RING
Consortium~\cite{Pellerin2017}. Clearly, the adoption of this technique calls
for numerical schemes able to operate on general polytopal elements.  This method
when applied to Cartesian grids has the advantage of maintaining a structured
grid away from the fracture network, where the sparsity of the linear system may
ease its numerical solution, but it does not allow local refinements (unless by
using hanging nodes, which increases computational complexity). In general it is a valid alternative to a direct triangulation provided the numerical scheme be robust with respect to the presence of small or high
aspect ratio elements.

We outline a possible algorithm for the case of a Cartesian background grid,
adopted in this work.  We start by creating a Cartesian mesh of rectangular
elements and compute the intersections among the edges of the grid and the
segments describing the fractures. This step is rather straightforward for
Cartesian grids. The intersection points can be easily sorted according to a
parametric coordinate to create the mesh of each fracture. Then, each cell cut
by one or more fractures is split into two, three or four polygonal sub cells as
follows: \textit{i)} for each point, the signed distance from the fracture is
computed, and \textit{ii)} points on the same side of the fracture are grouped,
and sorted in counter-clockwise direction.

To avoid non-convex cells the cells containing the fracture tips are split in
three by connecting the tip with the nodes of the edge that is crossed by the
prolongation of the fracture. However, in principle it would be possible to
consider a single cell with two coincident faces.

\subsection{Agglomeration}\label{subsec:agglomeration}

Polytopal grids can be generated by agglomerating simplicial elements
produced, for instance, by a constrained Delaunay procedure. For example,
in~\cite{Botella2015}, tetrahedra are agglomerated (and nodes moved) to try to
produce hexahedral elements in large part of the domain, with a twofold objective: on the one hand the reduction of the total number of degrees of freedom
and consequent reduction of computational complexity, on the other hand, the generation of a grid more suitable for finite volume schemes based on two-point flux
approximation (TPFA).

In a more general setting, agglomeration may join together elements whose value
of physical parameters are similar, with the final objective of reducing
computational cost, as well as eliminating excessively small elements. The
numerical method, however, should be able to operate properly on the possible irregularly
shaped and non-convex elements generated by the procedure. The technique is
clearly a post-processing one, since it requires to have a mesh to start with.
Its basic implementation is however rather simple and is similar to that used in
some multigrid solvers, like in~\cite{Jones2001a}.

In our case, PorePy has the capability to agglomerate cells based on two
different criteria: (\textit{i)} by volume, meaning that cells with small volumes
are grouped with neighbouring cells. This procedure continues until the new
created cells have volumes that are comparable with an averaged volume in the
grid. This procedure can be effective in presence of uniform physical data in
different part of the computational domain and in particular in presence of
fracture networks. In the case of highly variable data, e.g. permeability, the previous
procedure may not be effective since cells with very different properties may be
merged together. For this reason PorePy implements another strategy, (\textit{ii)}
based on the agglomeration in the algebraic multigrid method. It adopts a
measure of the strength of connections between DOFs to select the cells to be joined, based on a two-point flux approximation discretization, for
more details see \cite{Trottenberg2001,Fumagalli2016a}.
Examples of these strategies are given in
\cite{Fumagalli2016a,Fumagalli2017a,Fumagalli2017,Scialo2017,Nordbotten2018}.

\subsection{Voronoi} Voronoi grids are of particular interest for methods such
as Finite Volumes with TPFA, since they guarantee that the line connecting the
centroids of neighbouring cells is always orthogonal to the shared face. Under
this assumption the two point approximation of the flux is consistent if the
permeability tensor is diagonal. However, producing Voronoi diagrams that honour the internal interfaces represented by the fracture is not an easy
task, particularly for complex 3D configurations. An attempt in that direction has
been performed in~\cite{Pellerin2014} and~\cite{Berge2019a}.

\begin{figure}[htb]
    \centering
    \includegraphics[width=0.8\textwidth]{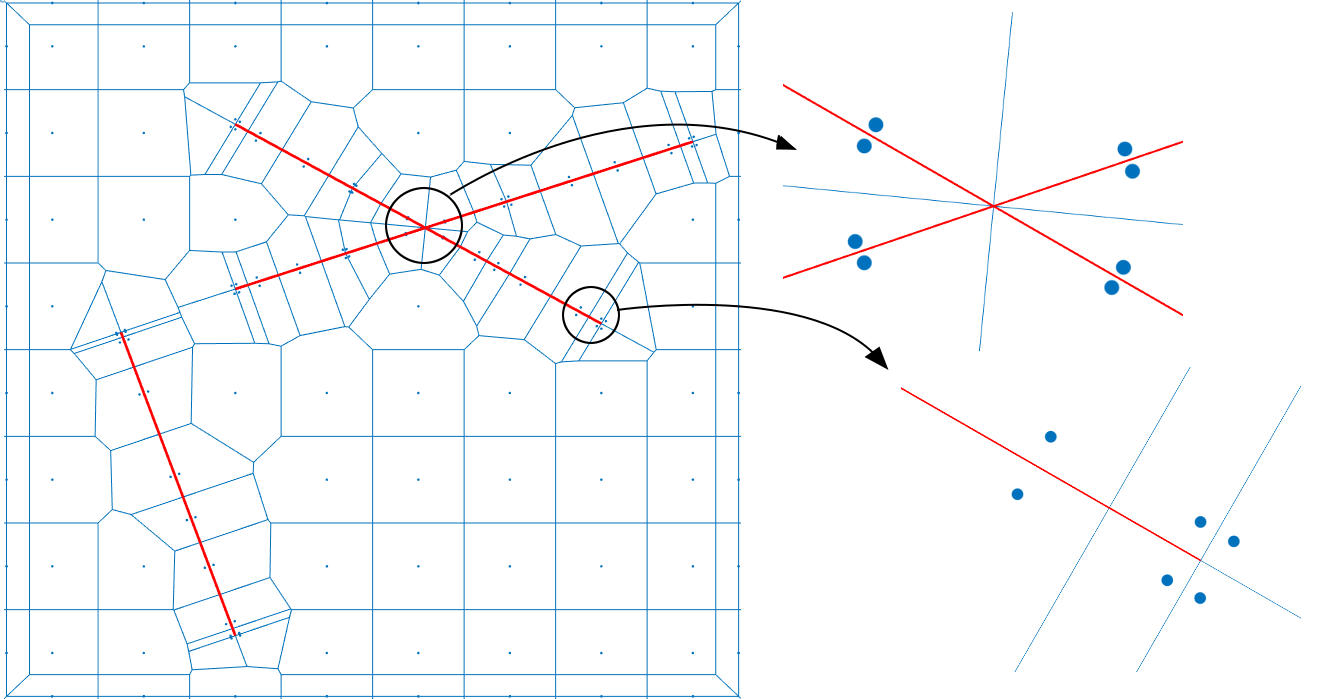}
    \caption{On the left, graphical representation of Voronoi grid with
    fractures. On the right, details on the construction for fracture
    intersection and fracture tip.}%
    \label{voro}
\end{figure}

In this work, limited to 2D cases, we generate Voronoi diagrams that honour the
geometry of the fractures and the boundaries of the domain by first creating a
Cartesian grid (see Section \ref{cart}) and positioning a seed at the centre of
cells not cut by the fractures. Then, we start from the discretization of the
fractures induced by the intersection with the background grid, and for each
fracture cell we position two seeds on opposite sides of the fracture at a small
distance $\delta$ as shown in Figure \ref{voro}. This will create a Voronoi cell
with a face exactly on the fracture. The same technique is used to obtain
boundary faces in the desired position. Close to each fracture tip
$\bm{x}_T$ we position four seeds in $\bm{x}_T\pm \delta_1 \bm{n}
\pm \delta_2 \bm{t}$ where $\bm{n}$ and $\bm{t}$ are the normal and
the tangent unit vectors to the fracture and $\delta_{1,2}$ are user defined
distances. This ensures that the fracture is honoured up to the tip and has the
correct length. Similar strategies are applied at fracture intersections. The position of the seeds and faces close to the intersections
is also shown in Figure \ref{voro}. Note that with this strategy the Voronoi cells far from fractures are rather regular, since they reflect the underlying Cartesian grid.

An advantage of Voronoi grids is that faces are planar and cells are convex by
construction. However, an important drawback is that the number of faces per
cell can be quite large. Moreover, as pointed out before, the construction of a
constrained grid in general realistic configurations is an open problem.

\section{Numerical results}\label{sec:numerical_results}

In this section, we present two test cases to show the performances and
the potentiality of the previously introduced algorithms. In particular, in the
first test case we have a setting where the permeability experiences a high
variation between neighbouring cells. In the second test case a network of
fractures is considered with different types of intersections: in this case the
challenge is more related to the geometrical complexity to create the
computational grid.
In both test cases, coarsening procedures are used to reduce the computational
cost of the simulations.

\subsection{Heterogeneous porous medium: layers from SPE10}\label{subsec:spe10}

\begin{figure}[htb]
    \centering
    \subfloat[L4 original ]
    {\includegraphics[width=0.33\textwidth]{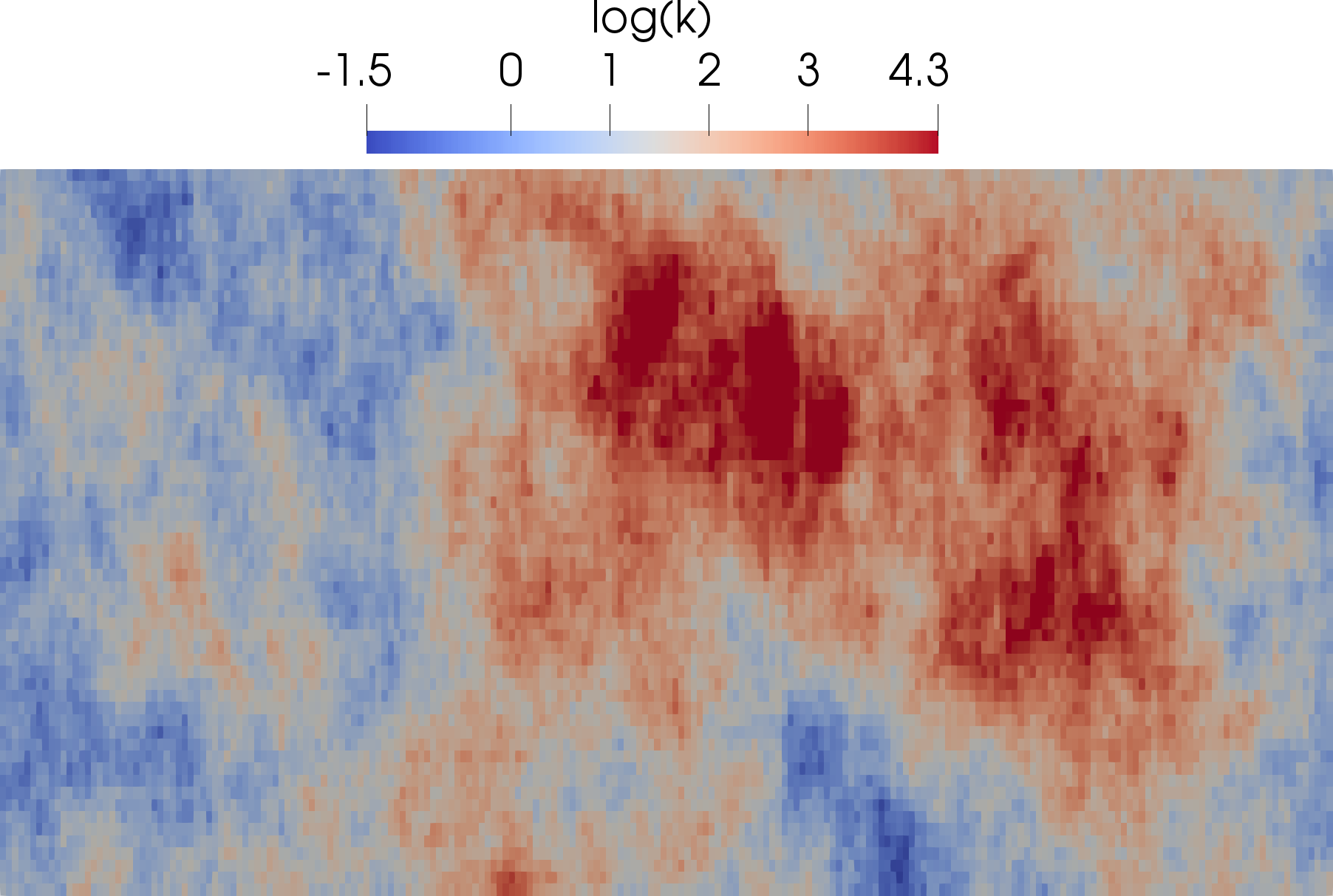}}%
    \hfill%
    \subfloat[L4 coarsened with arith. mean]
    {\includegraphics[width=0.33\textwidth]{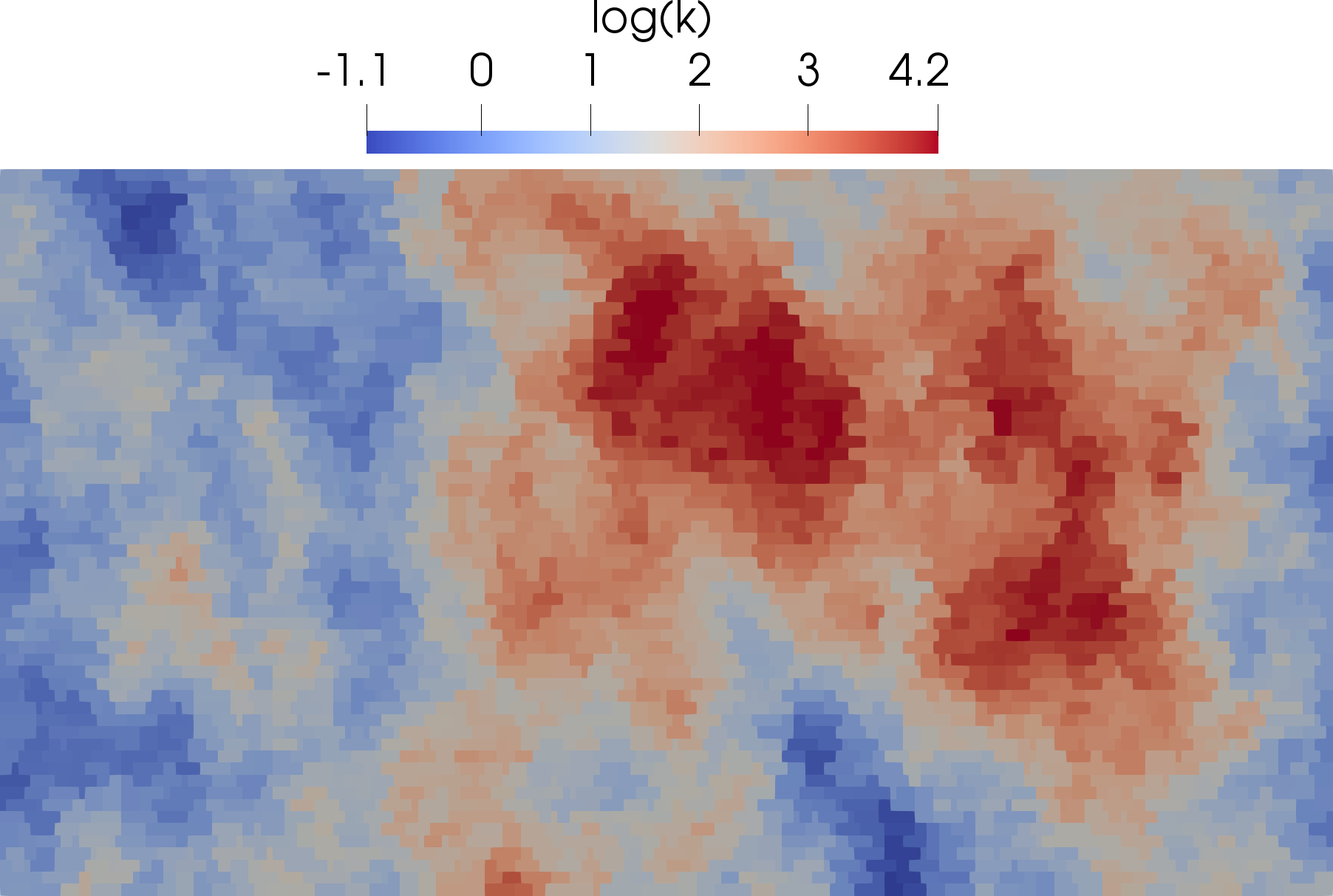}}%
    \hfill%
    \subfloat[L4 coarsened with harm. mean]
    {\includegraphics[width=0.33\textwidth]{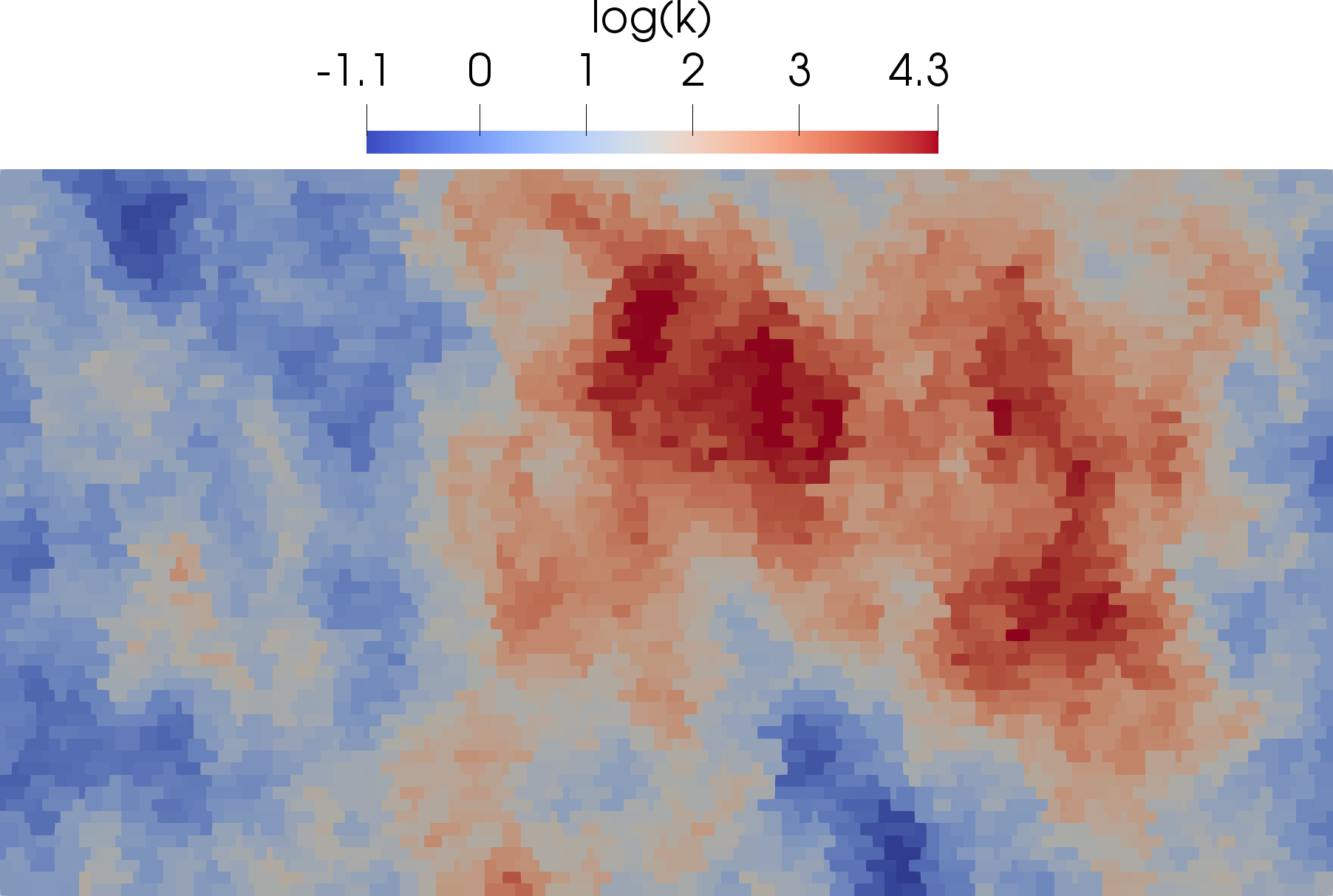}}\\
    \subfloat[L35 original ]
    {\includegraphics[width=0.33\textwidth]{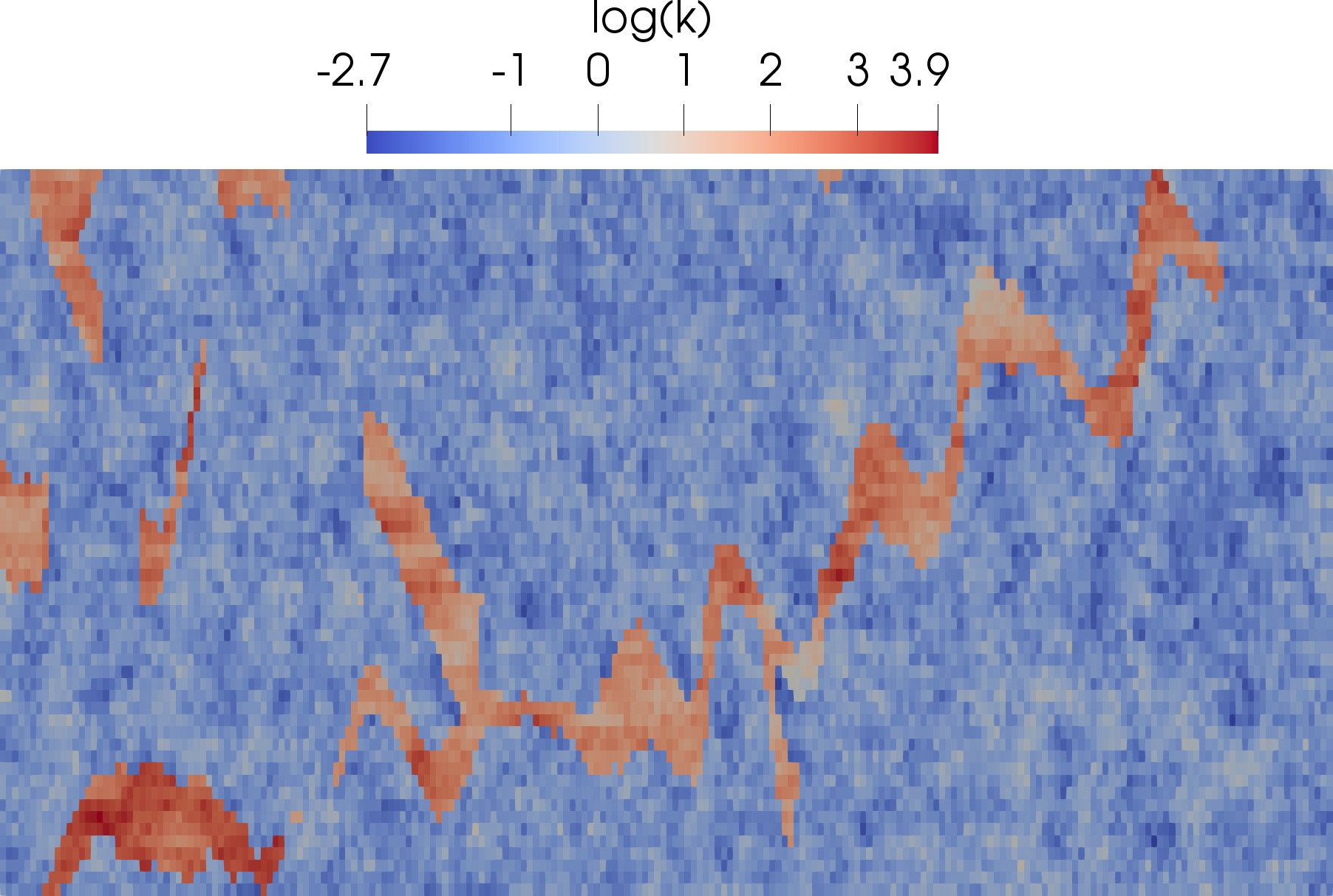}}%
    \hfill%
    \subfloat[L35 coarsened with arith. mean]
    {\includegraphics[width=0.33\textwidth]{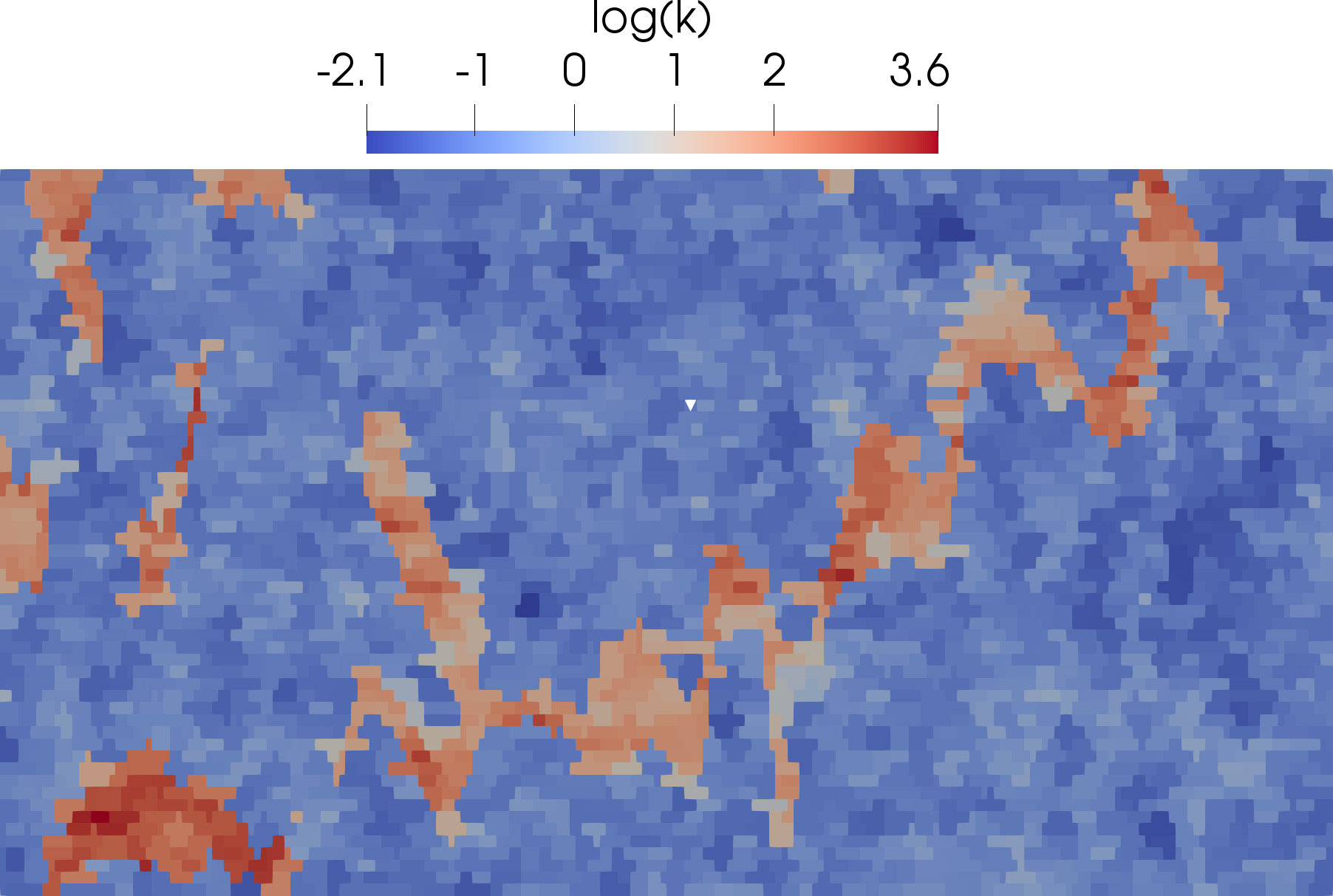}}%
    \hfill%
    \subfloat[L35 coarsened with harm. mean]
    {\includegraphics[width=0.33\textwidth]{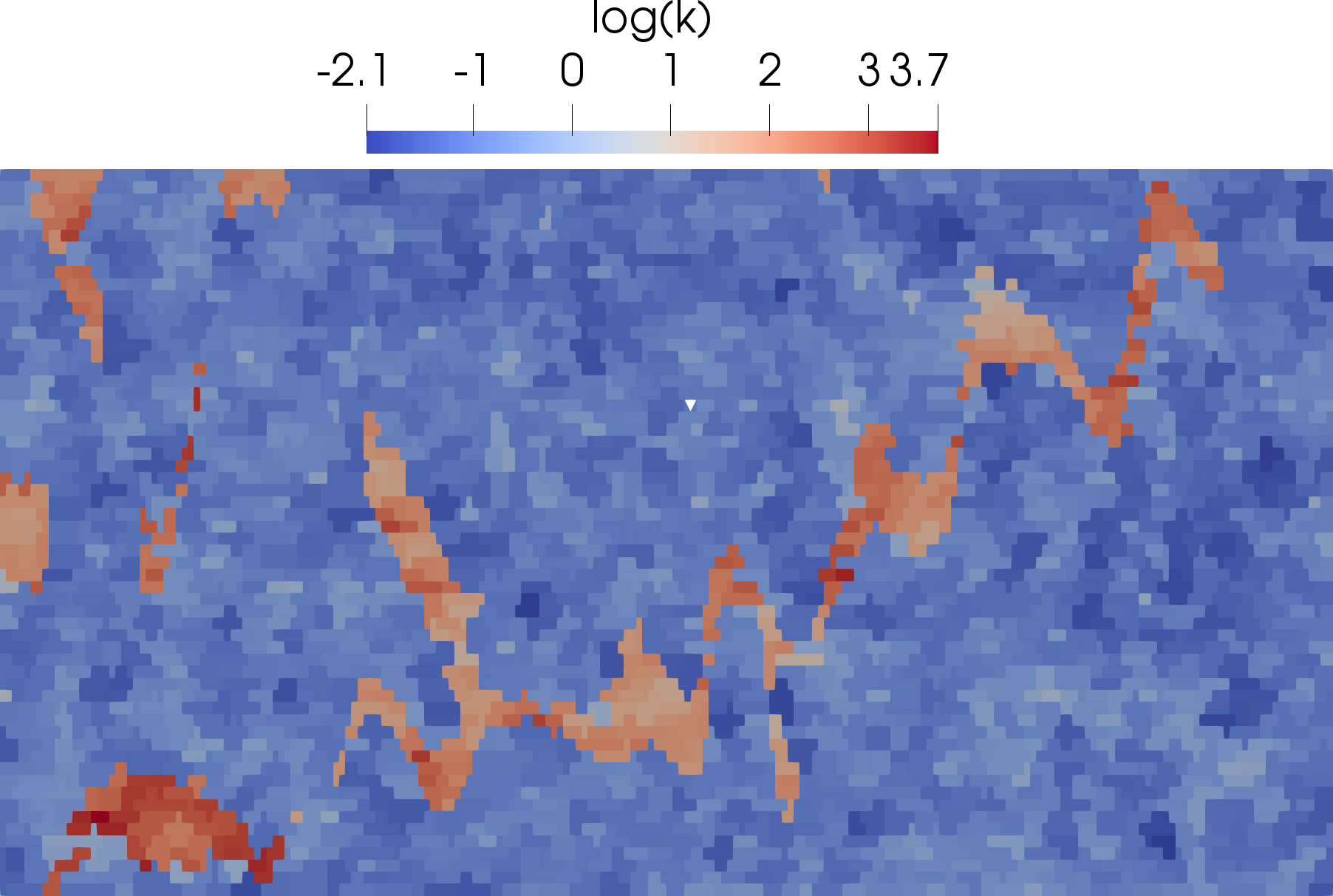}}%
    \caption{Permeability field for the test case of Subsection
    \ref{subsec:spe10} for layer 4 on the top and 35 on the bottom.
    On the left the reference values, on the
    centre and right the values obtained after the clustering with arithmetic
    and harmonic mean, respectively. The values are given in $\log_{10}$.}
    \label{fig:spe10_perm}
\end{figure}

The aim of this test case is to validate the effectiveness of the MVEM scheme in
presence of highly heterogeneous permeability. We consider two
distinct layers of the SPE10 \cite{Christie2001} benchmark problem, in particular layer 4
and 35 (by starting the numeration from 1), from now on denoted as L4 and L35,
respectively. The main difference between them is
that the latter has distinctive channels of high permeability which are not
present in layer 4. The permeability is assumed to be scalar in each cell, and
each layer is composed by a computational grid of
$60\times220$. Figure \ref{fig:spe10_perm} on the
left shows the permeability fields for both layers. Note that in both cases permeability spans about six order of magnitude.

To lighten the computational effort, we apply a coarsening
procedure to group cells and obtain a smaller problem. Starting from square
cells the algorithm creates cells by considering the procedure in Subsection
\ref{subsec:agglomeration} and, for each coarsened cell, the associated permeability will be computed
in two different ways: as the arithmetic and harmonic average. The former is more
suited for flow parallel to layers of different permeability, while for orthogonal flow
the harmonic average gives more realistic results. For a more detailed
discussion see \cite{Nordbotten2011}. We consider both approaches, see Figure
\ref{fig:spe10_perm} on centre and right, which represents the coarsened
permeability of both layers by considering the arithmetic and harmonic means.
For layer 4 the figures look similar, while for layer 35 the channels for the
coarsened problem with harmonic mean are narrower than the original ones and
than those obtained in the coarsened grid with arithmetic mean.

In Table \ref{tab:aree_facceSPE}, we summarize the geometric properties of the
grids obtained by means of cells clustering for the two layers. We can observe
that the area of the cells and the average number of faces per cell is similar
in the two cases, however, in layer 35 we have slightly more elongated elements
on average, reflecting the channelized permeability field. The aspect ratio is
estimated using the area of the cells, the maximum distance between points is
rescaled so that square cells (or equilateral triangles, see Section
\ref{subsec:benchmark}) have aspect ratio 1.

\begin{table}\centering
 \begin{tabular}{|l|c|c|c|c|c|c|c|c|c|}
 \hline
  \multirow{2}{*} &\multicolumn{3}{c|}{aspect ratio} & \multicolumn{3}{c|}{cell area} &\multicolumn{3}{c|}{$n_\mathrm{faces}$}\\     \cline{2-10}
    & average & min & max  & average & min & max & average & min & max\\ \hline
    L4 & 2.37& 1.50 & 4.37 & 108 & 37.2 &242&  12.2 &  6 & 20\\ \hline
    L35 & 2.37 & 1.13 & 5.83 & 111 & 37.2 & 297 &  12.2 &  6 & 22\\ \hline
 \end{tabular}
\caption{Average, minimum and maximum value of cell area and number of faces per cell for the six grids employed for Test case \ref{subsec:spe10}}\label{tab:aree_facceSPE}
\end{table}

We impose a pressure gradient from left to right with synthetic values 1
and 0, respectively. The other boundaries are sealed with homogeneous Neumann
conditions.

To compare the accuracy of the proposed clustering techniques, we compute the
errors in the pressure with respect to the problem on the original grid solved
with a two-point flux approximation scheme \cite{Aavatsmark2007a,Friis2009},
which, in this case since the grid is $K$-orthogonal, is consistent and
converges quadratically to the exact solution, thus can be considered as a valid
reference. We name this solution ``reference'' and we indicate the pressure as
$p_{\rm ref}$. The error is computed as
\begin{gather*}
    err = \frac{\norm{\Pi_{\rm ref} p - p_{\rm ref}}_{L^2(\Omega)}}{\norm{p_{\rm ref}}_{L^2(\Omega)}}
\end{gather*}
where $\Pi_{\rm ref}$ is the piecewise constant projection operator that maps
from the current grid to the reference one. Due to the clustering procedure its
construction is rather straightforward, since the cells of the original mesh are nested in
the coarsened one.
We can notice that the errors obtained for the layer 4 with both averaging
procedure are comparable and around $4\%$, which can be acceptable in most of
real applications. In the case of layer 35 the situation is more involved, in
fact the arithmetic mean gives an error of approximately $3.5\%$ while the
harmonic mean of $13\%$. We can explain this big discrepancy by noticing that,
when a channel of high permeability is composed
by few cells in its normal direction, during the coarsening procedure it is
possible that some of these cells are grouped with the surrounding lower
permeability cells. The harmonic mean will bring the permeability value of the coarsened cell closer to the lower value than the higher, dramatically changing the connectivity properties of the obtained permeability field. This can be noticed in the permeability field reported in Figure \ref{fig:spe10_perm}, suggesting that harmonic averaging can be unsuited for parallel flow in strongly channeled domains.
\begin{figure}
    \centering
    \subfloat[L4 original ]
    {\includegraphics[width=0.33\textwidth]{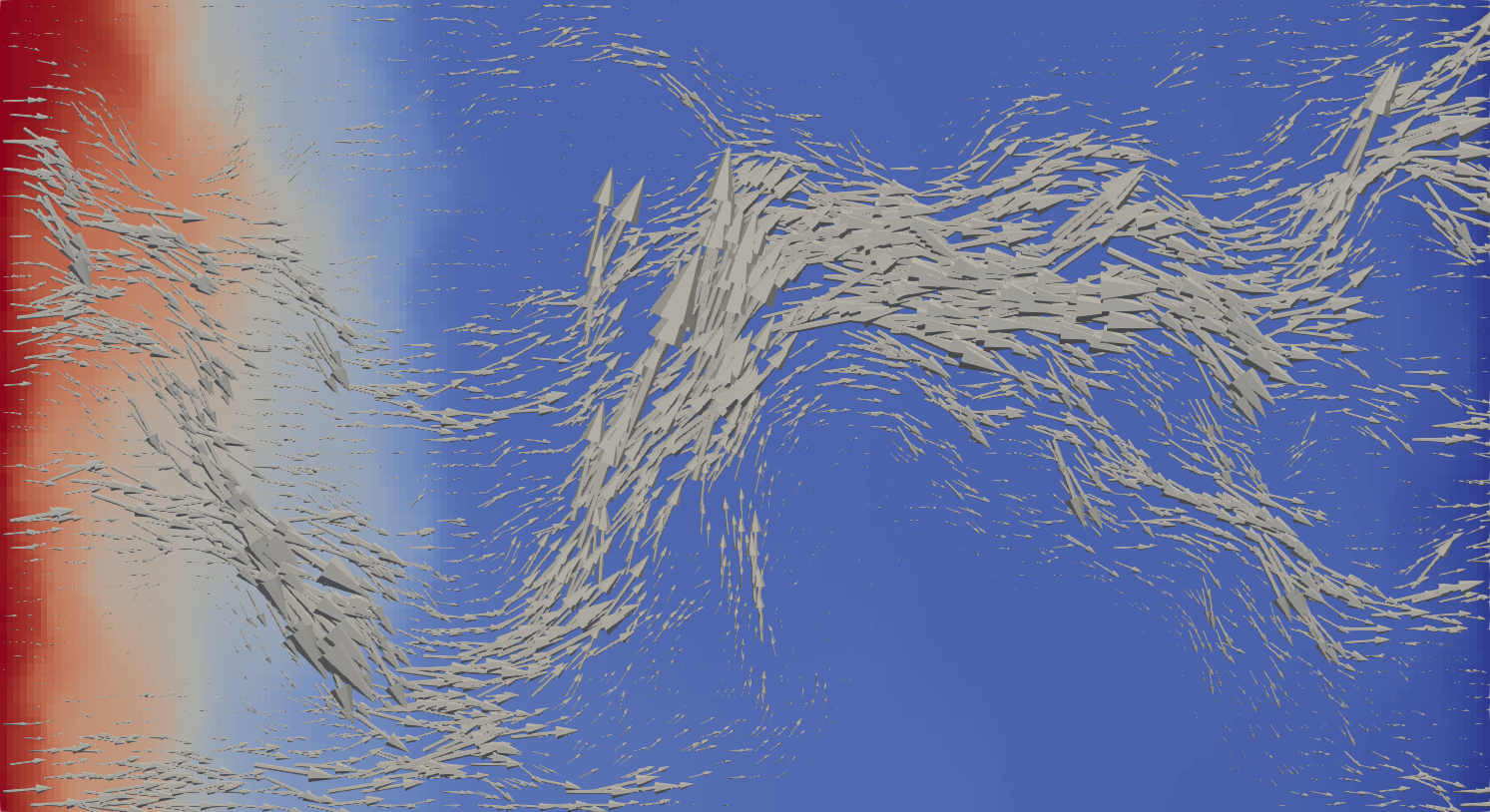}}%
    \hfill
    \subfloat[L4 coarsened with arith. mean]
    {\includegraphics[width=0.33\textwidth]{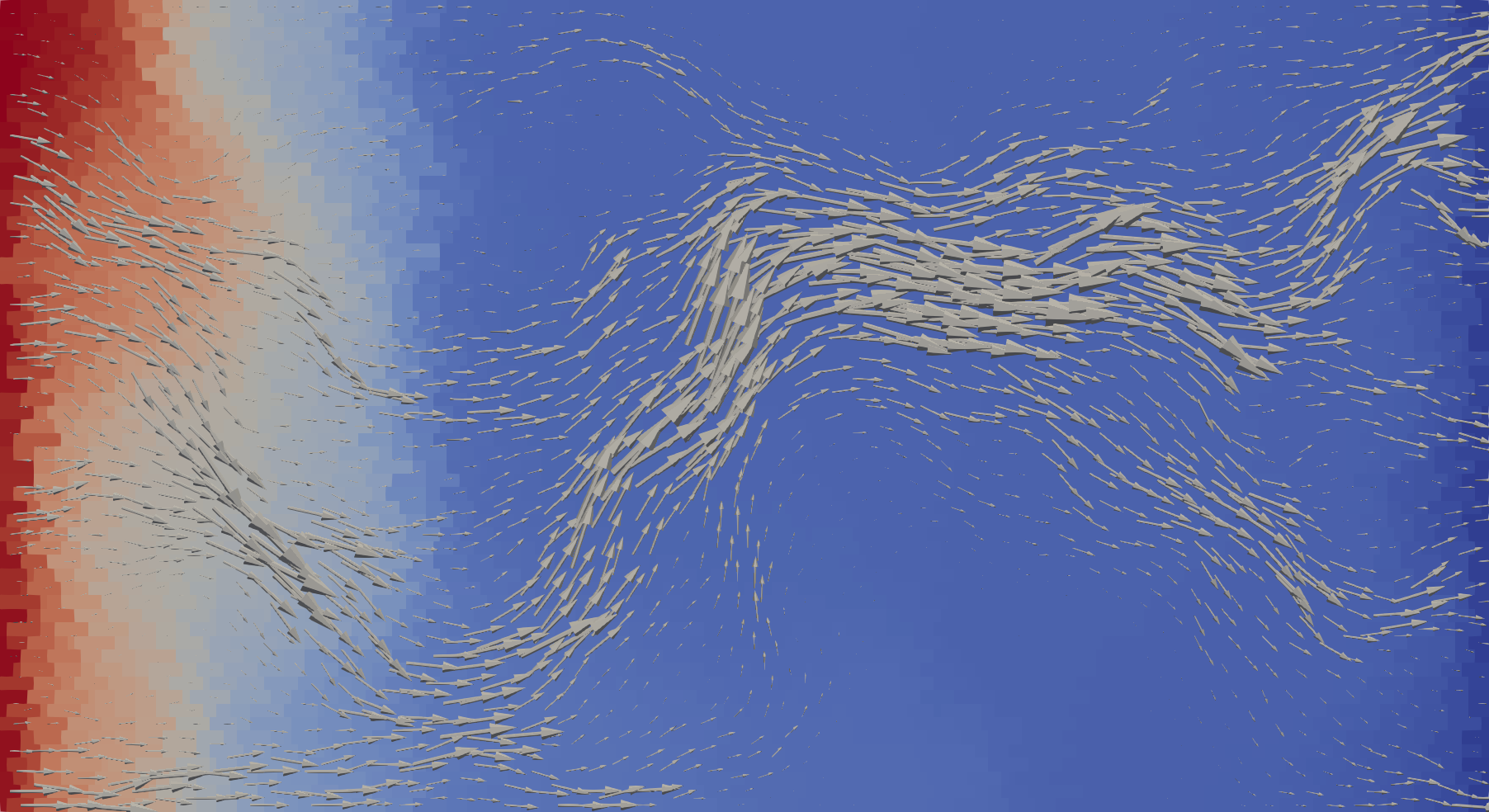}}%
    \hfill
    \subfloat[L4 coarsened with harm. mean]
    {\includegraphics[width=0.33\textwidth]{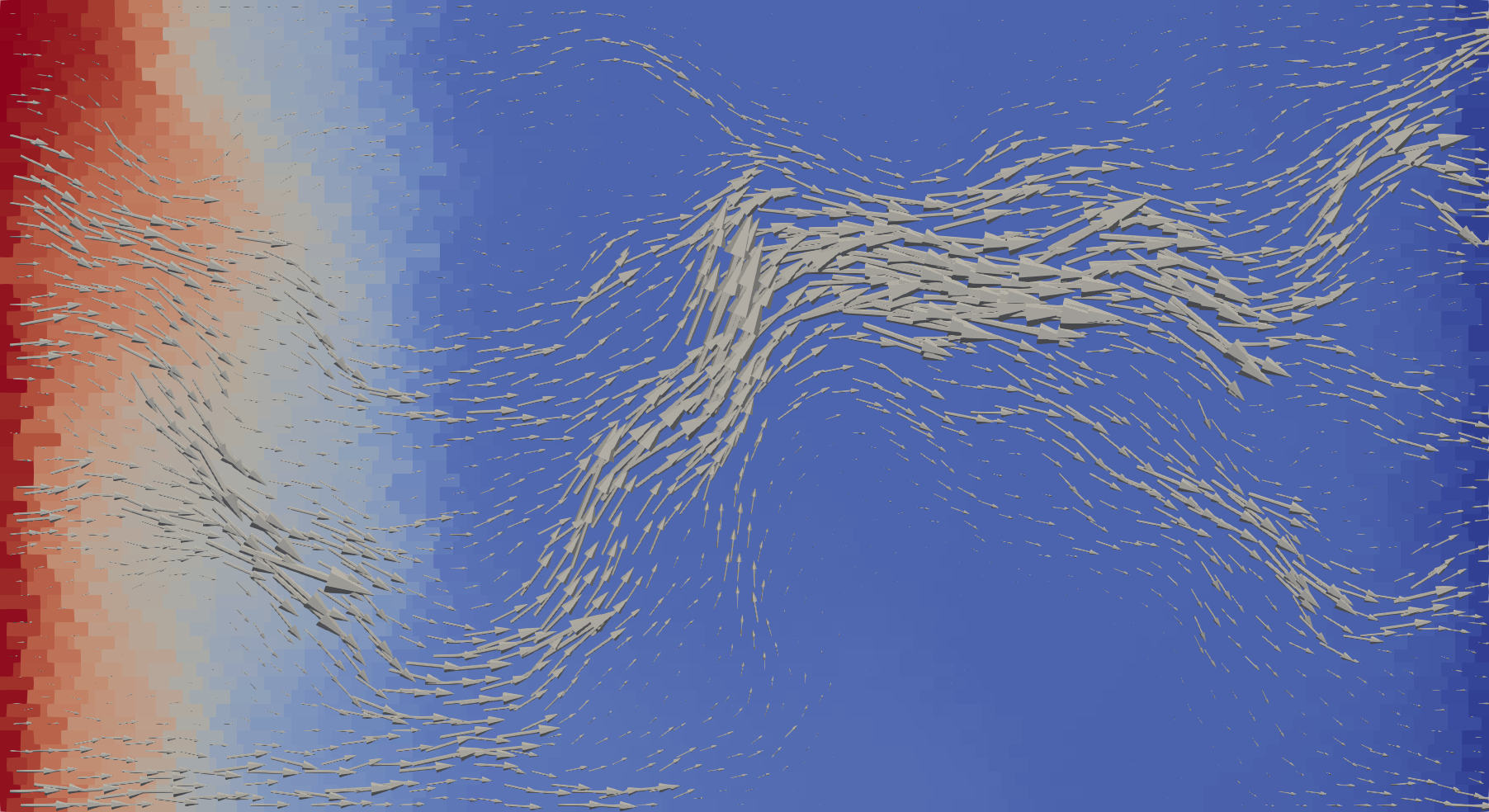}}\\
    \subfloat[L35 original ]
    {\includegraphics[width=0.33\textwidth]{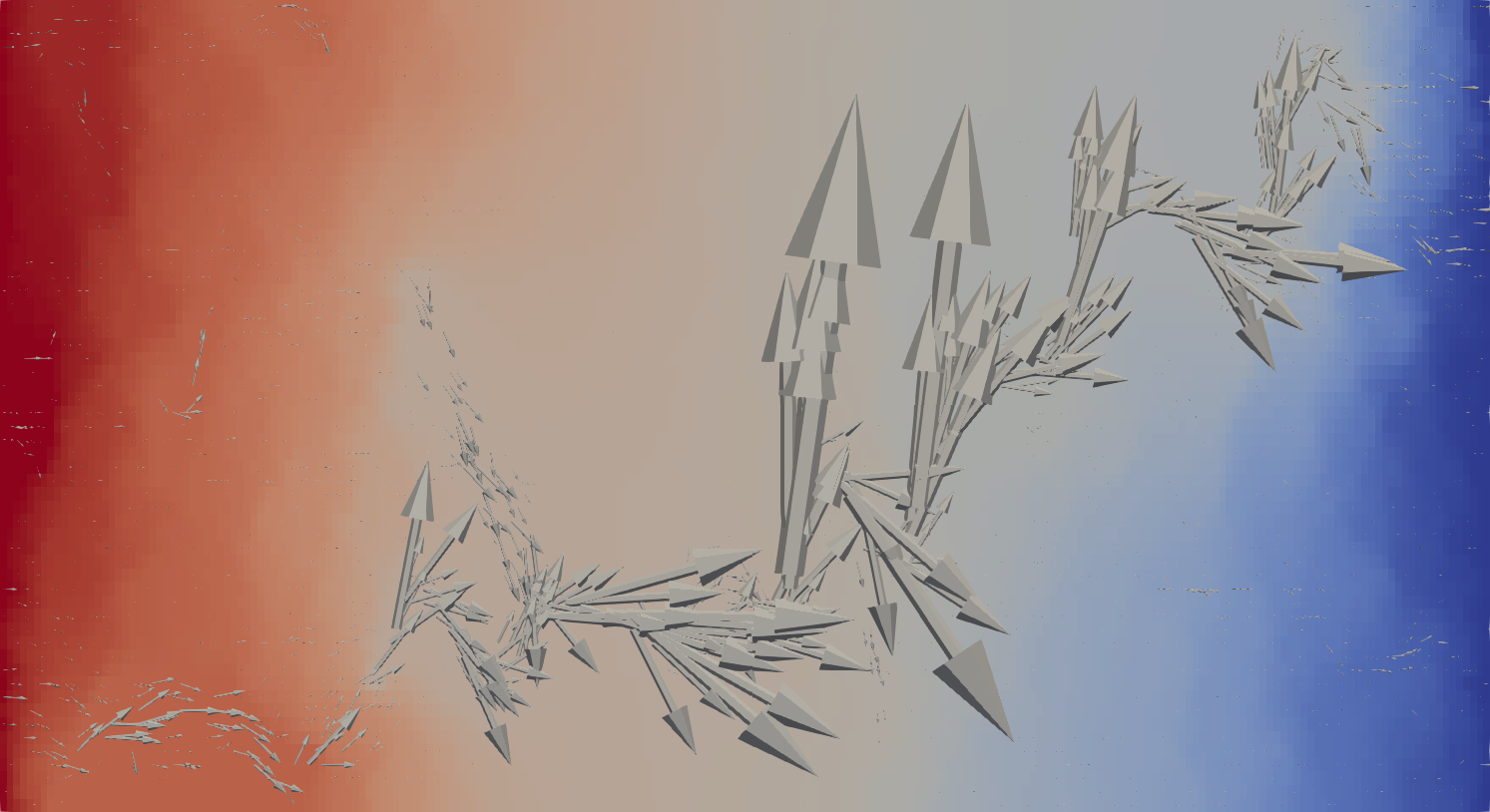}}%
    \hfill
    \subfloat[L35 coarsened with arith. mean]
    {\includegraphics[width=0.33\textwidth]{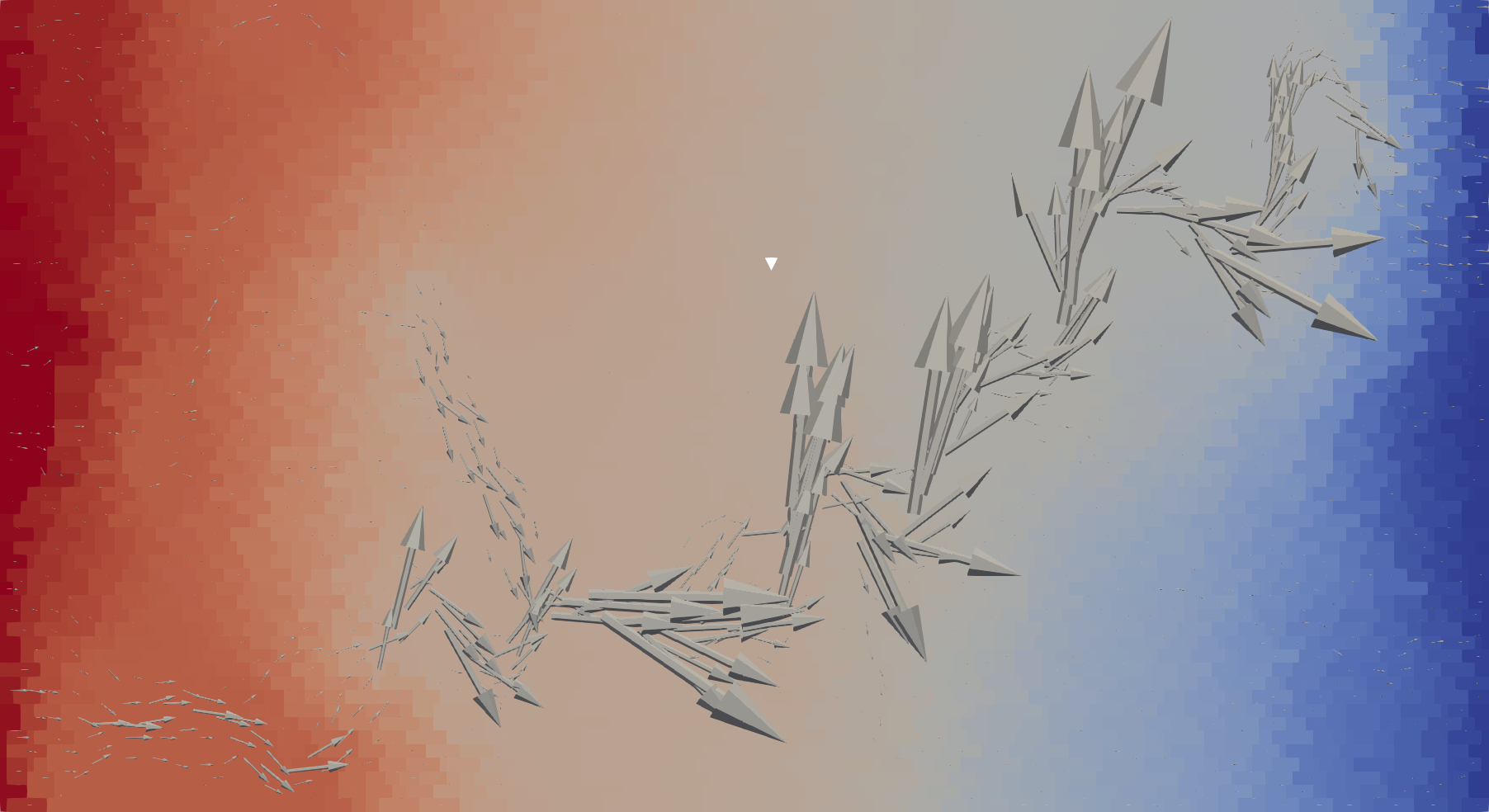}}%
    \hfill
    \subfloat[L35 coarsened with harm. mean]
    {\includegraphics[width=0.33\textwidth]{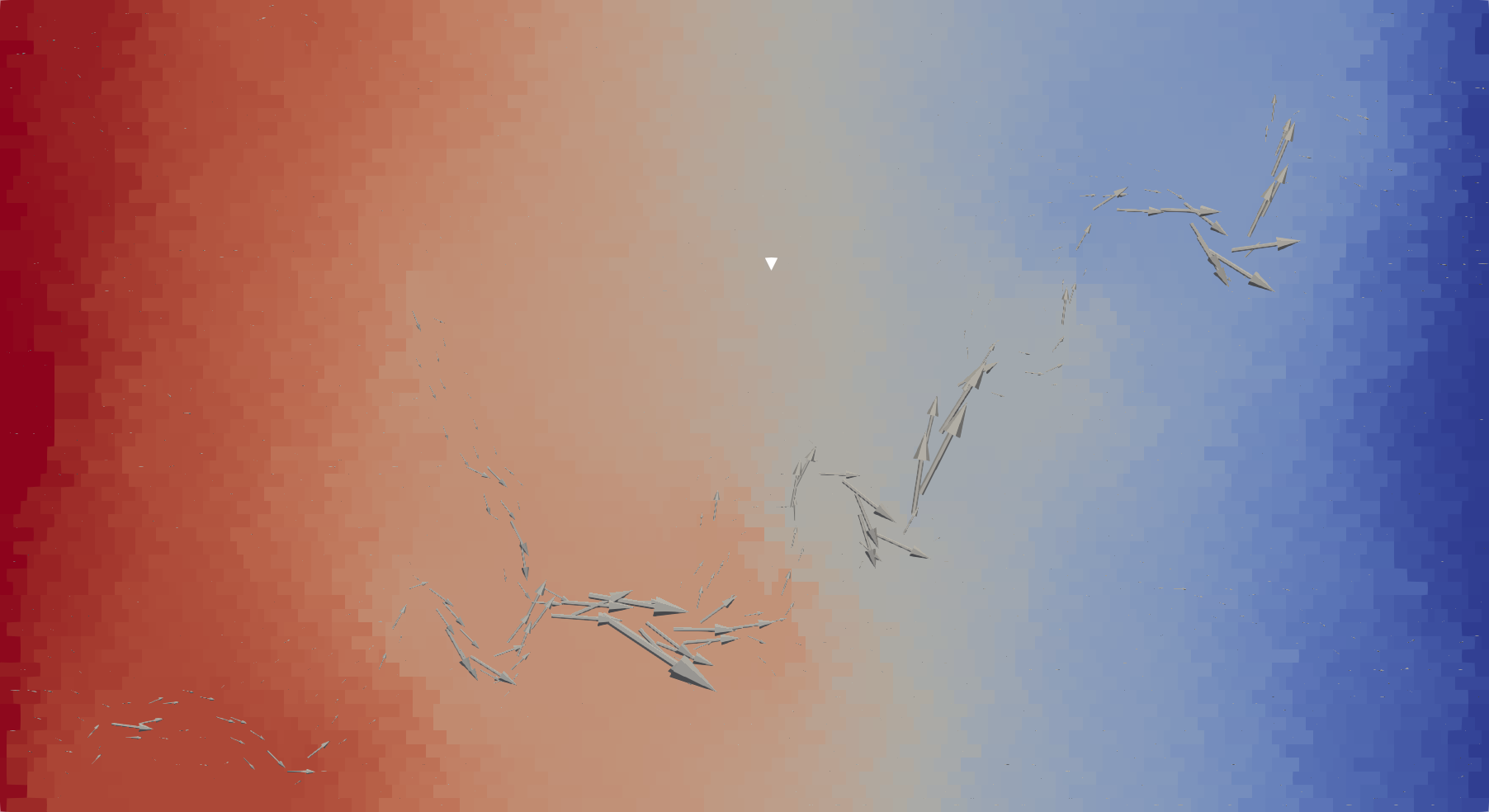}}%
    \caption{Pressure and Darcy velocity fields for the test case of Subsection
    \ref{subsec:spe10} for layer 4 on the top and 36 on the bottom. On the left the reference solution, on the
    center and right the values obtained after the clustering with arithmetic
    and harmonic mean, respectively.
    The arrows are scaled by the same value in each layer and the pressure ranges from 0 to 1,
    blue to red respectively.}
    \label{fig:spe10_pressure_35}
\end{figure}

Figure \ref{fig:spe10_pressure_35} shows the pressure fields for both layers and
for the two approaches. On top of the pressure fields the Darcy
velocity is also represented with grey arrows. We notice that for layer 4
pressures and velocities look very similar, while for layer 35 the
pressure field and velocity of the coarsened problem with harmonic mean look
quite different compared with the reference solutions as well as that obtained with the coarsened strategy that uses the arithmetic mean.

To improve the effectiveness of this approach, a local numerical upscaling technique could be
considered to compute a more representative value of the permeability for
grouped cells. However, in this case we might expect a higher computational
cost. See \cite{Durlofsky2003} for a more detailed presentation of upscaling techniques.

To conclude this test case, let us now analyze the properties of the system matrix to verify what is the
impact of element size and shape in the different cases. Note that the problem is in mixed form and our analysis considers the entire saddle point matrix. Since the number of
unknowns is not exactly the same after grid agglomeration we describe matrix
sparsity by means of the average number of non-zero entries per row $\overline{n}$, computed as
\begin{gather*}
    \overline{n}=\frac{n_z}{N_\mathrm{DOF}},
\end{gather*}
where $n_z$ is the number of non-zero entries and $N_\mathrm{DOF}$ is the number
of unknowns. Moreover we will compare the values of condition number $K(A)$
estimated by the method \texttt{condest} provided by
\texttt{Matlab}\textsuperscript\textregistered. In Table \ref{tab:matrici_caso1}
we consider the two layers, L4 and L35, and by ``mean K'', ``harmonic K'' we
identify the averaging of permeability in the coarse cells, the arithmetic and
harmonic mean respectively. This choice has no impact on the matrix size or
sparsity but may result in different condition numbers. We can observe that the
four matrices are very similar in terms of size, sparsity and condition number, and that the large number of faces per element reflects in the average number of entries per row.
It can be also observed that mesh coarsening is slightly more effective in layer
L35 due to its channelized permeability distribution.
\begin{table}\centering
 \begin{tabular}{|l|c|c|c|c|c|}
 \hline
    & $N_\mathrm{DOF}$ & $N_\mathrm{cells}$ & $N_\mathrm{faces}$ & $\overline{n}$ & $K(A)$ \\ \hline
    L4 (mean K)     & 16345 & 2269 & 14076 & 22.17 &  8.29e+06 \\ \hline
    L4 (harmonic K) & 16345 & 2269 & 14076 & 22.17 &  8.44e+06 \\ \hline
    L35 (mean K)    & 16010 & 2210 & 13800 & 22.53 &  8.29e+06 \\ \hline
    L35 (harmonic K)& 16010 & 2210 & 13800 & 22.53 &  8.39e+06 \\ \hline
 \end{tabular}
\caption{Matrix properties for Test case \ref{subsec:spe10}}\label{tab:matrici_caso1}
\end{table}

\subsection{Fracture network}\label{subsec:benchmark}

\begin{figure}[tbp]
    \centering
    \includegraphics[width=0.315\textwidth]{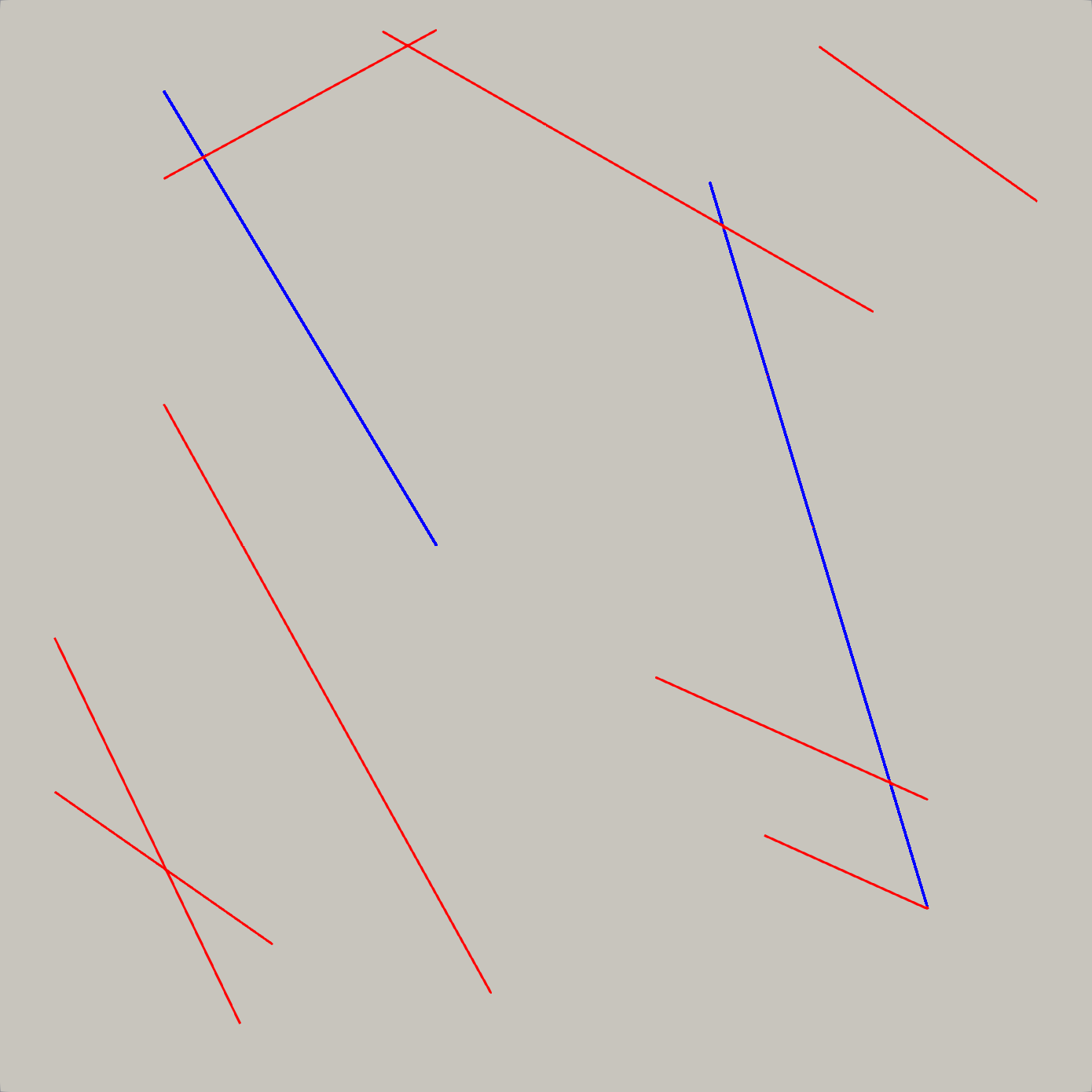}
    \caption{Geometry of the domain for the benchmark used in Subsection \ref{subsec:benchmark}.}
    \label{fig:anna_domain}
\end{figure}

This test case considers the Benchmark 3 of the study \cite{Flemisch2016a}
presented in Subsection 4.3. Our objective is to study the impact of the grid on
the solution quality provided by the MVEM. The domain contains a fracture
network made of 10 fractures and 6 intersections, one of which is of $L$-shape.
For the detailed fracture geometry, we refer to the aforementioned work.
See Figure \ref{fig:anna_domain} for a representation of the problem geometry.

\begin{figure}[hbt]
\centering
\subfloat[Delaunay (1275, 276)]{
  \includegraphics[width=0.315\textwidth]{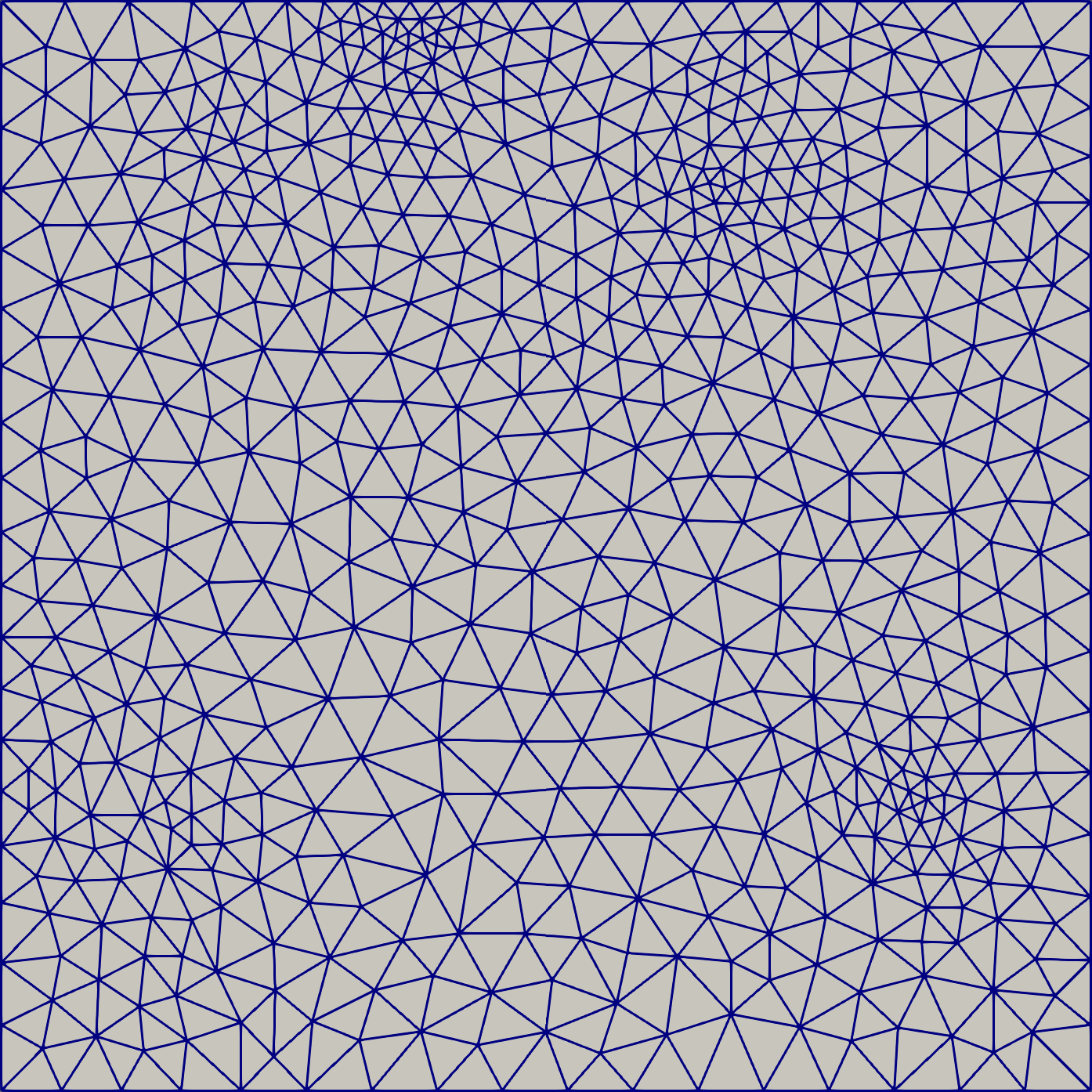}
}
\subfloat[Cartesian cut (1296, 579)]{
  \includegraphics[width=0.315\textwidth]{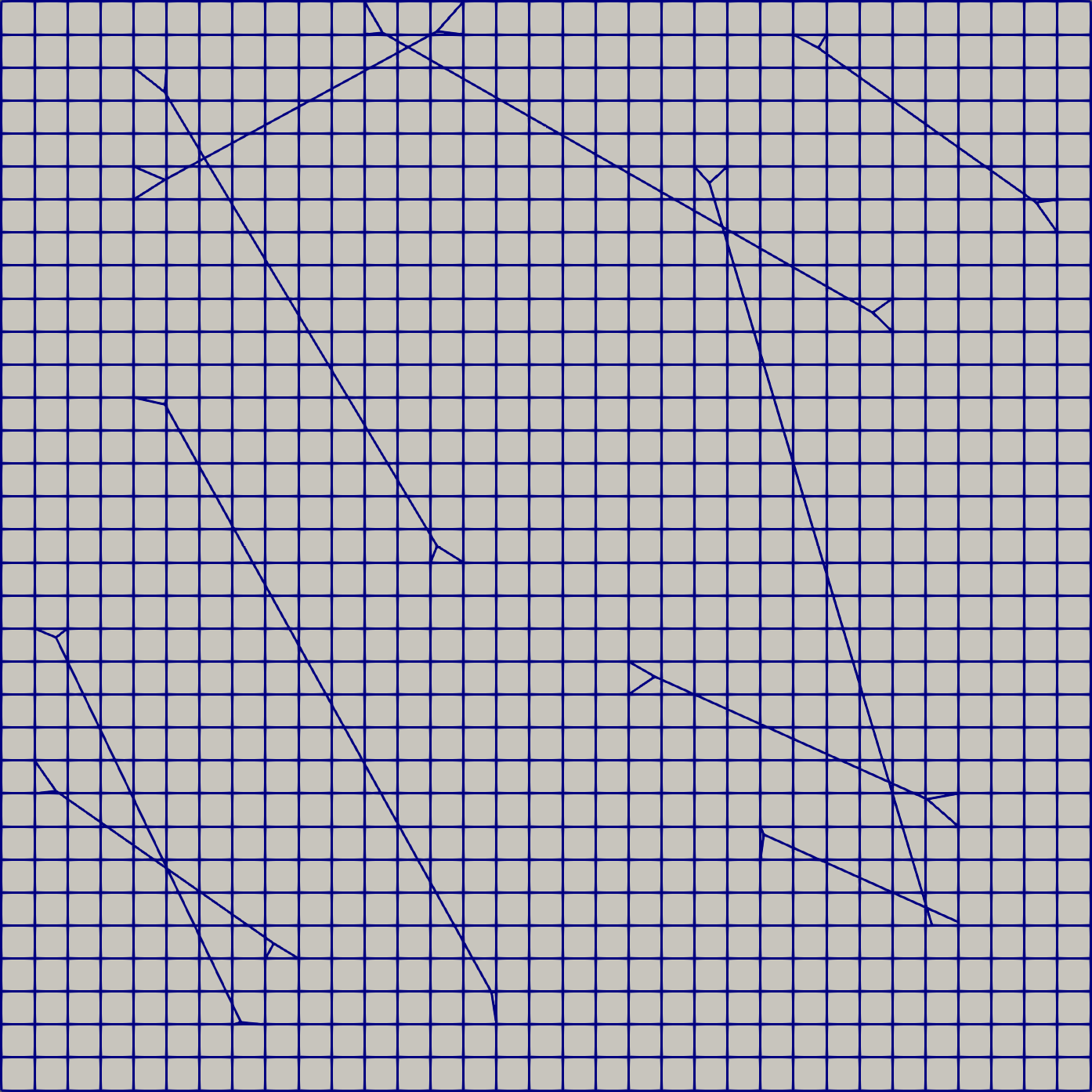}
}
\subfloat[Voronoi (1521, 633)]{
  \includegraphics[width=0.315\textwidth]{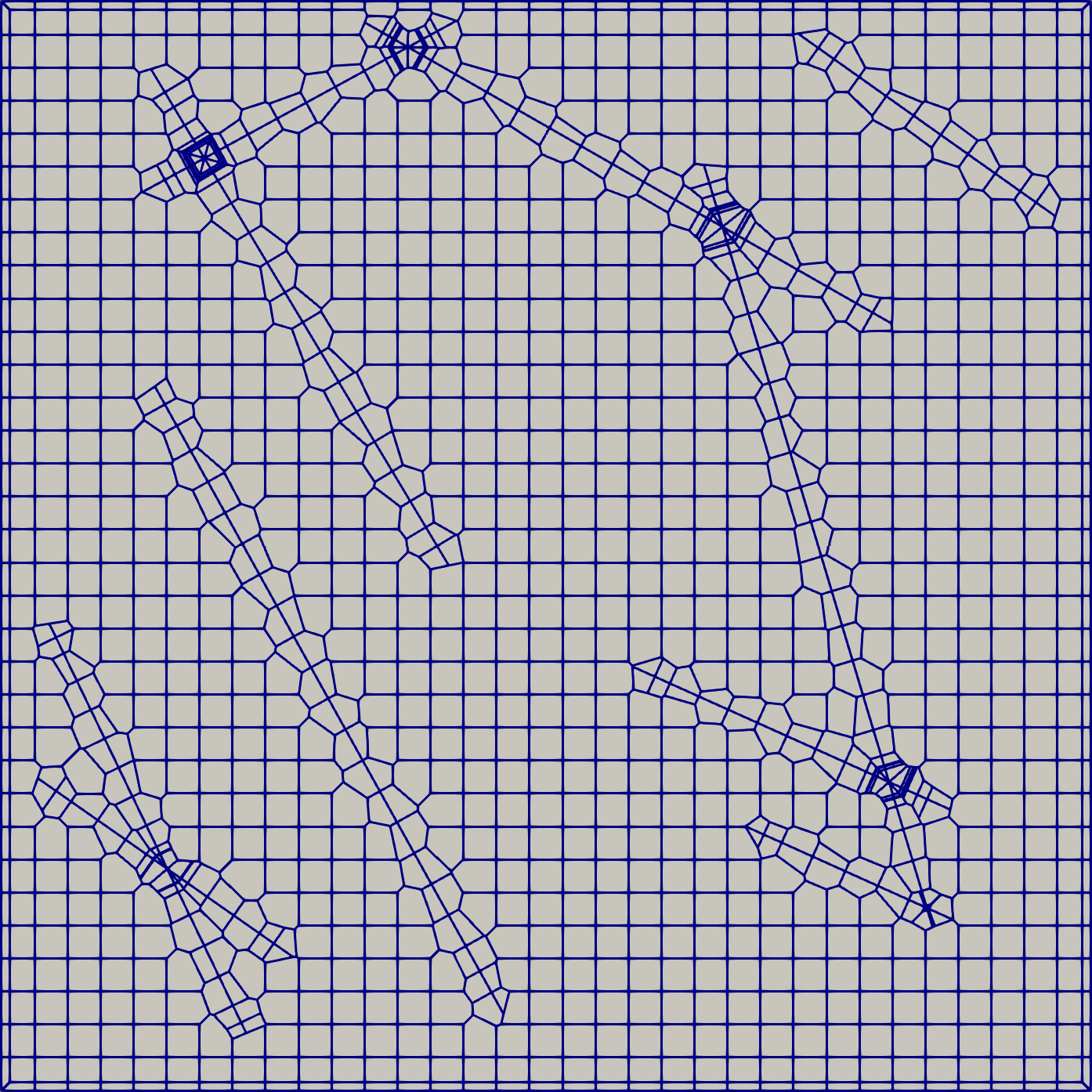}
}
\\
\subfloat[Delaunay coarse (1098, 276)\label{fig:case2_grid_delaunay_coarse}]{
  \includegraphics[width=0.315\textwidth]{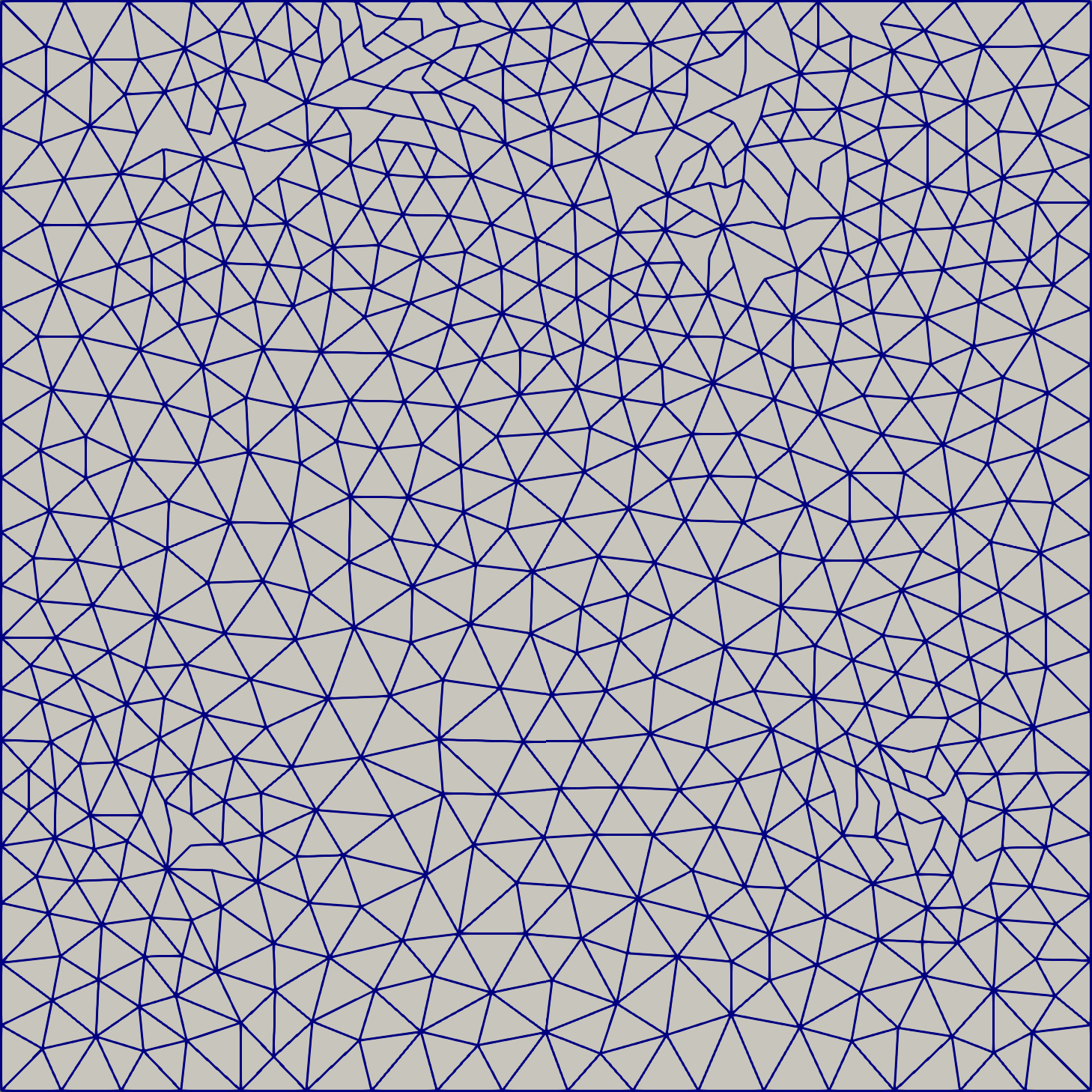}
}
\subfloat[Cartesian cut coarse (1053, 579)]{
  \includegraphics[width=0.315\textwidth]{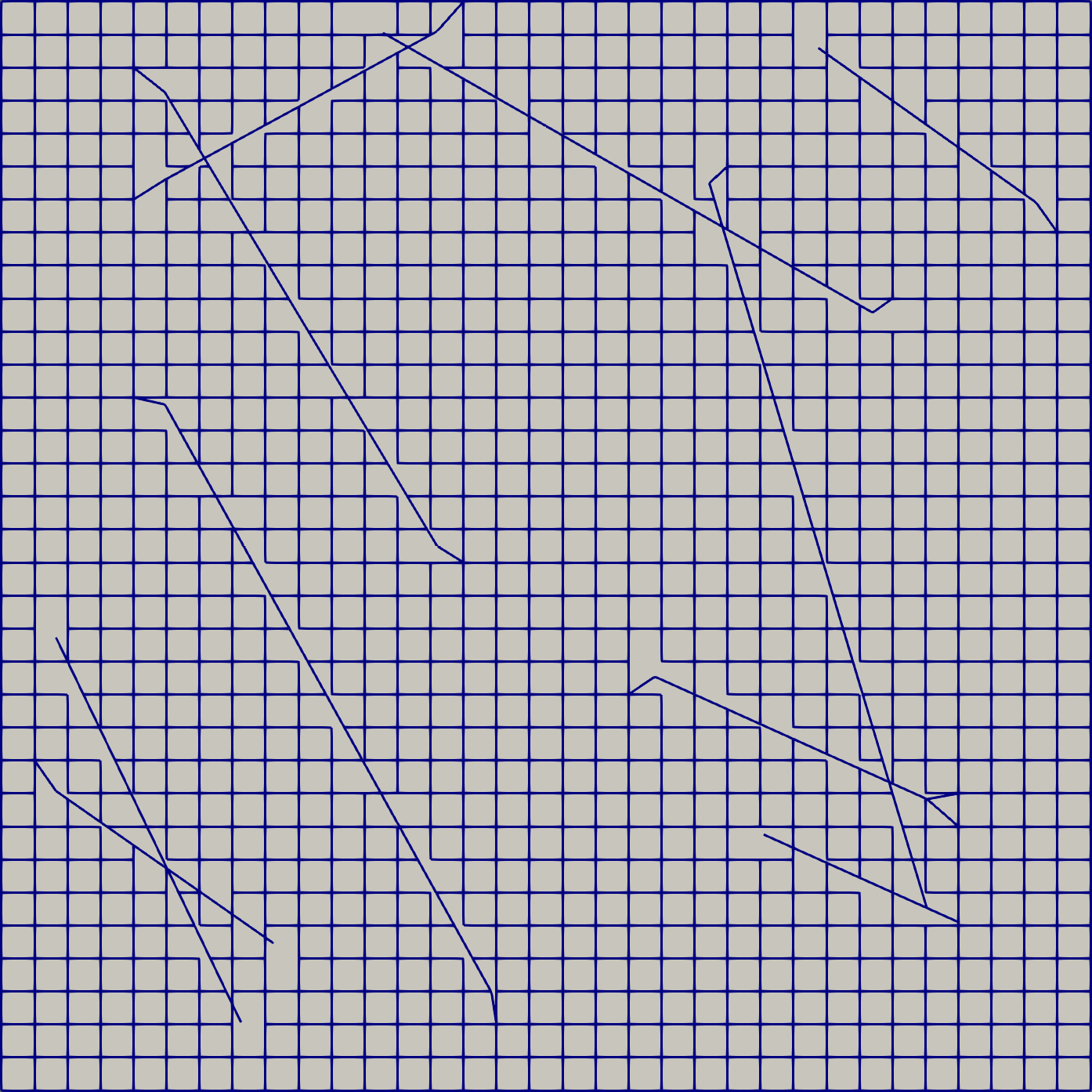}
}
\subfloat[Voronoi coarse (1109, 633)]{
  \includegraphics[width=0.315\textwidth]{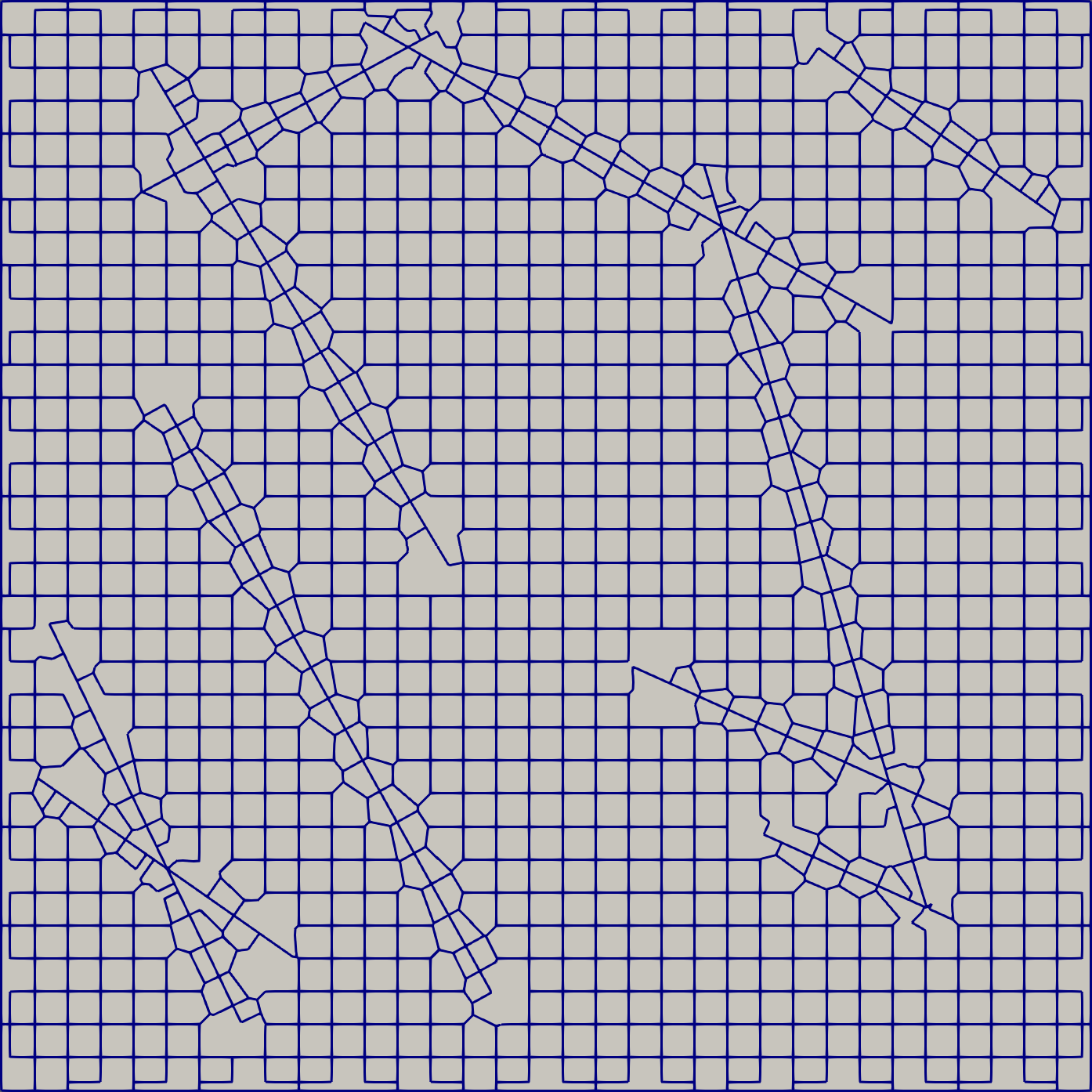}
}
\caption{Benchmark 3 of Subsection \ref{subsec:benchmark}: Fracture network on top left, on the others the grids for
different approaches. In the brackets the number of cells (bulk, fracture).}
\label{fig:anna_grids}
\end{figure}

We consider three types of grids:
Delaunay, Cartesian cut, and Voronoi. Since the fracture network may create
small cells, on top of these three grids a coarsening algorithm is
used to agglomerate cells of small volume. These cells are merged with neighbouring cells,
trying to obtain a more uniform cell size in the grid.
The Delaunay grid is created by the software Gmsh \cite{Geuzaine2009}, tuned to provide
high quality elements in proximity of small fracture branches or almost
intersecting fractures.
The six different grids we are considering are reported in Figure \ref{fig:anna_grids} along with the
number of cells associated to the rock matrix and fractures.
We see that some of the agglomerated elements have internal cuts, in particular
for Figure \ref{fig:case2_grid_delaunay_coarse}, and for all the clustered grids
we have cells that are not shape regular and in some cases not even star-shaped.
For classical finite elements or finite volumes we might expect low quality
results.

\begin{figure}[hbt]
\centering
\subfloat[Delaunay ]{
  \includegraphics[width=0.315\textwidth]{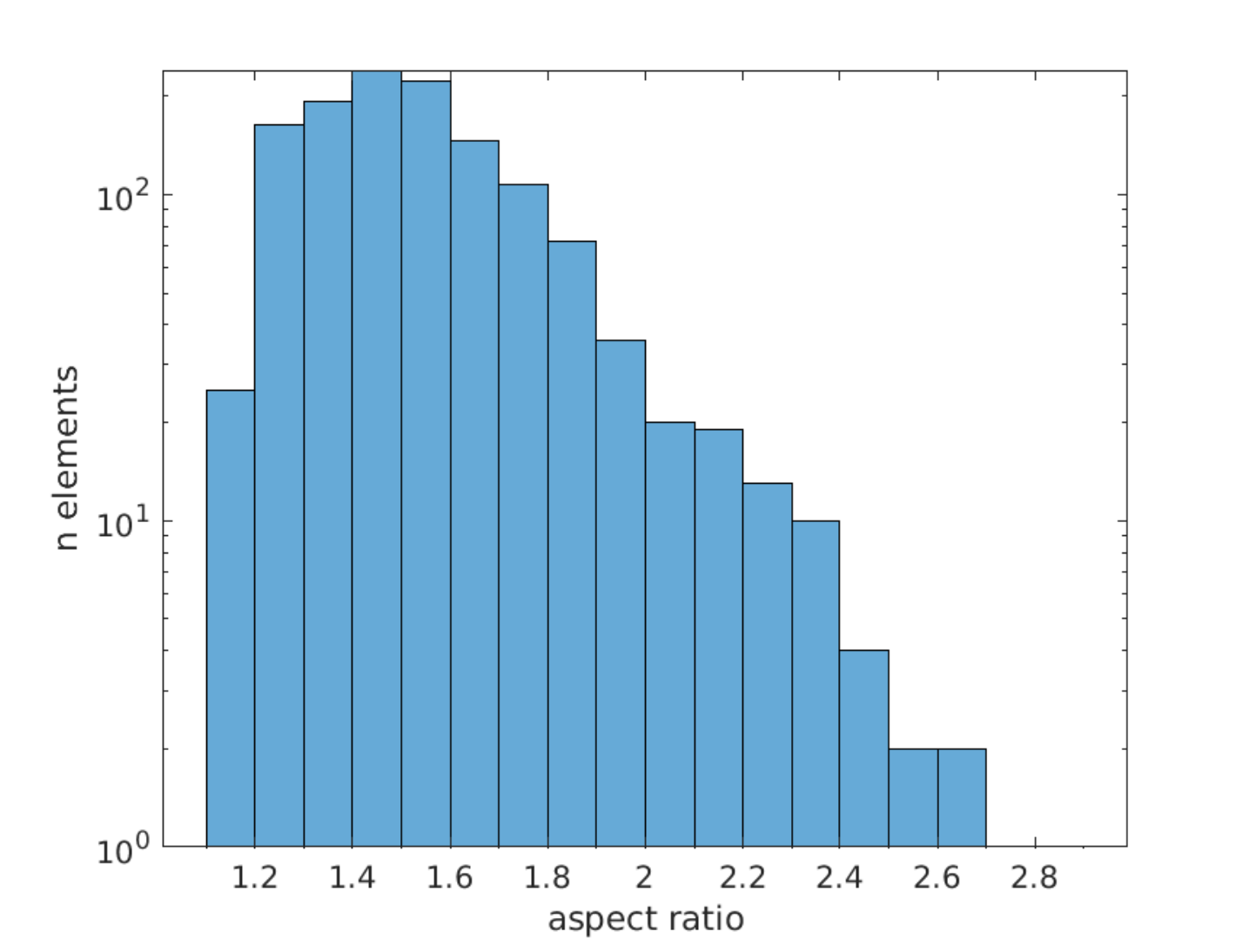}
}
\subfloat[Cartesian cut ]{
  \includegraphics[width=0.315\textwidth]{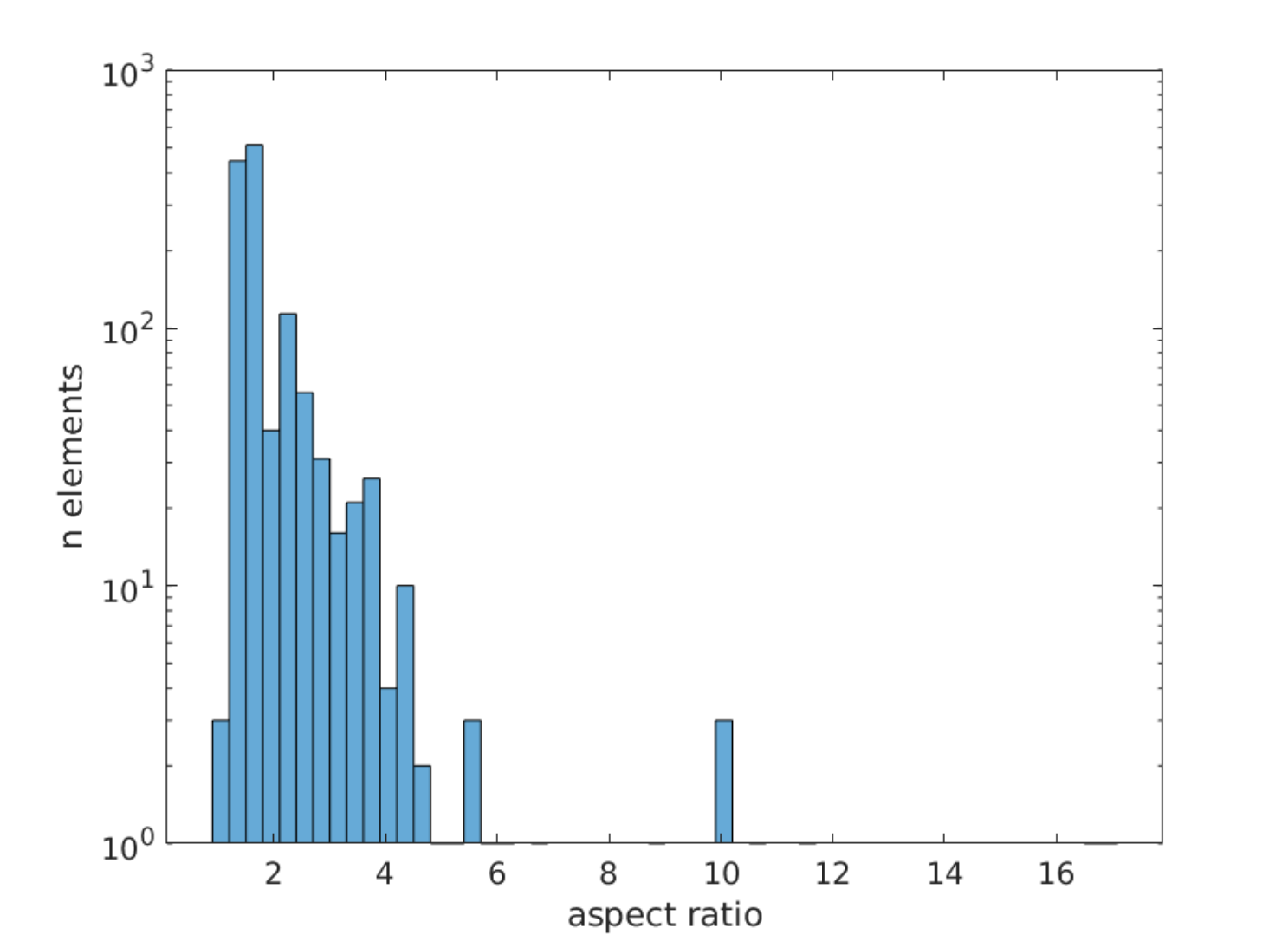}
}
\subfloat[Voronoi]{
  \includegraphics[width=0.315\textwidth]{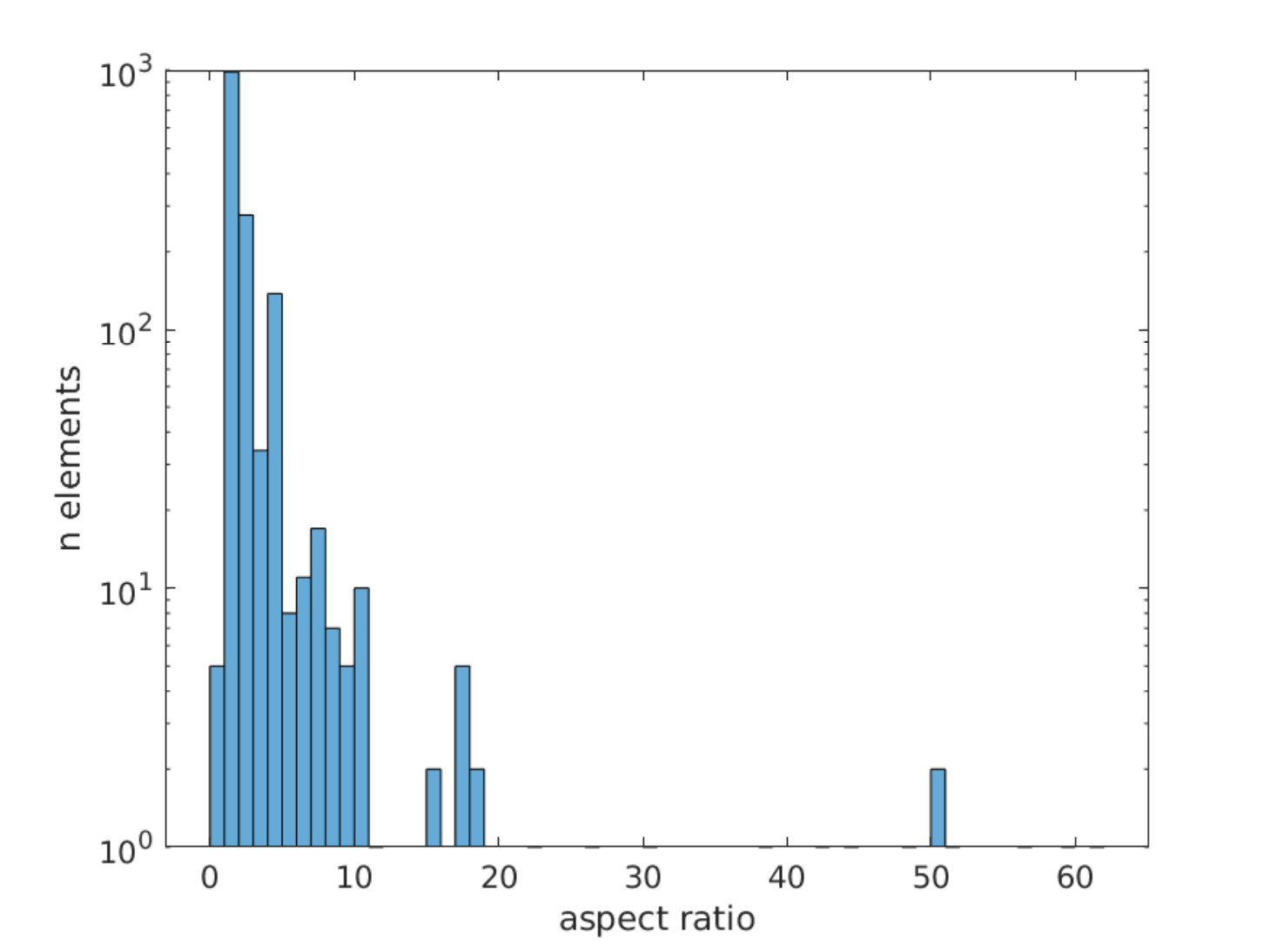}
}
\\
\subfloat[Delaunay coarse ]{
  \includegraphics[width=0.315\textwidth]{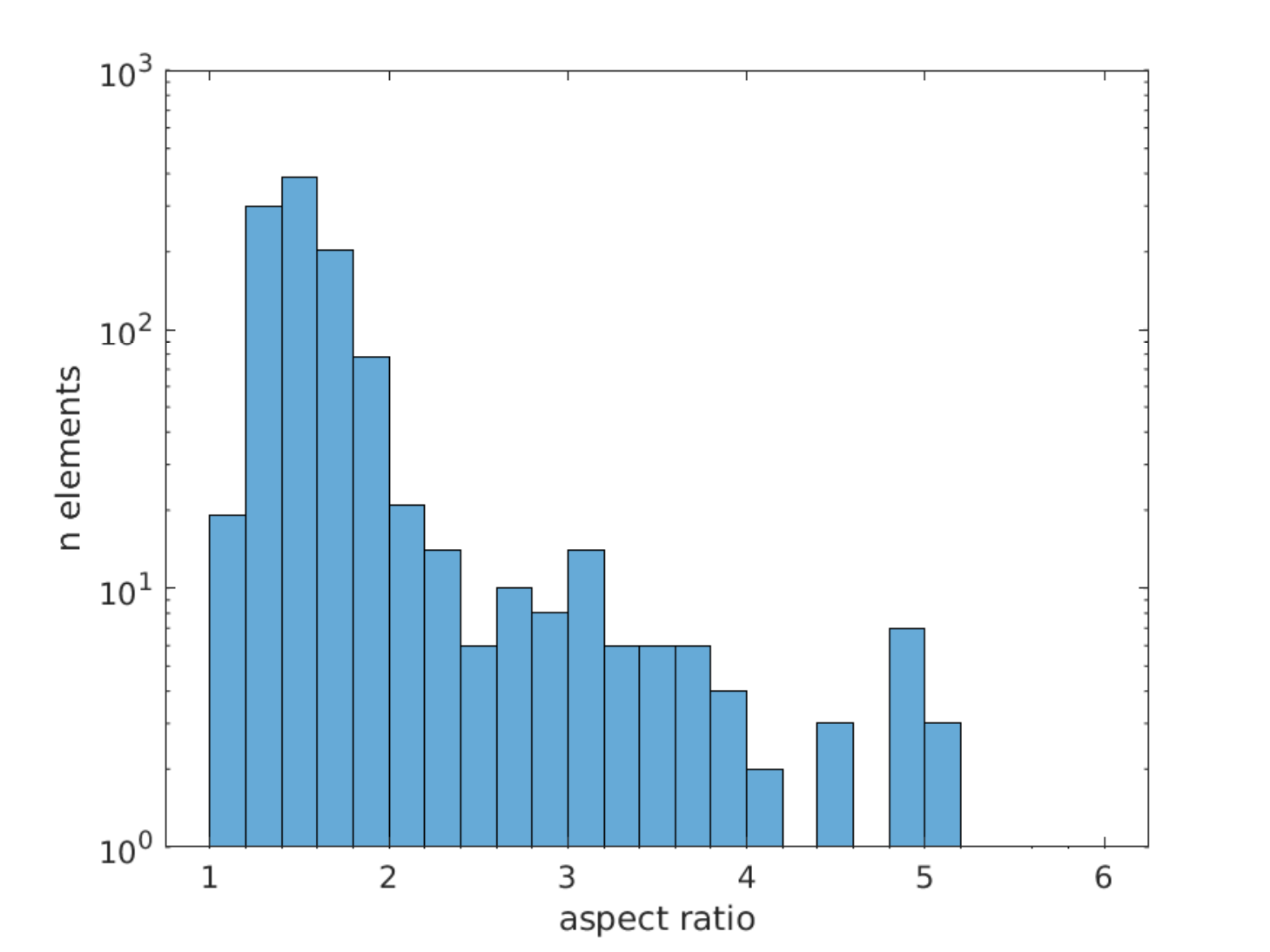}
}
\subfloat[Cartesian cut coarse ]{
  \includegraphics[width=0.315\textwidth]{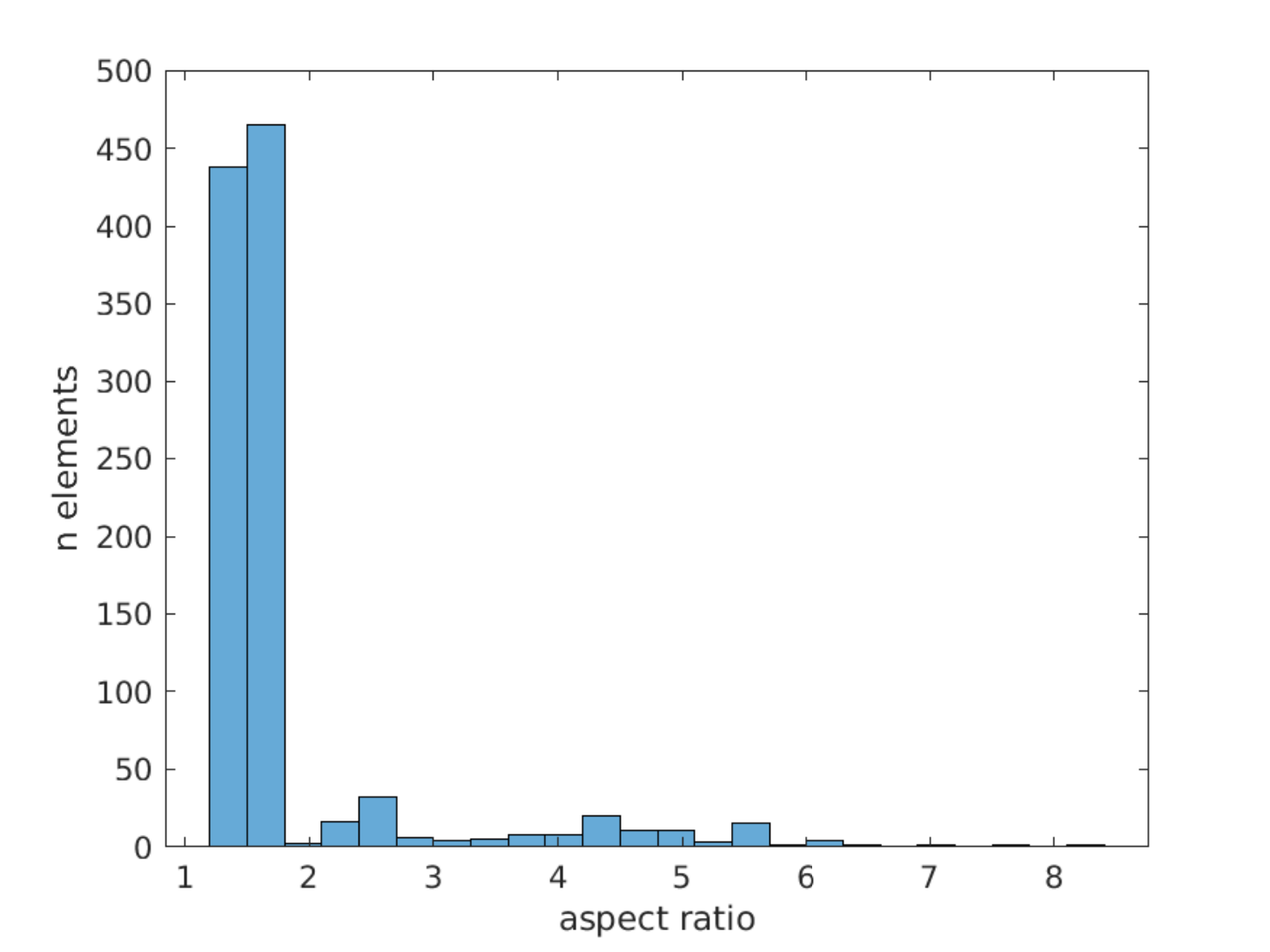}
}
\subfloat[Voronoi coarse ]{
  \includegraphics[width=0.315\textwidth]{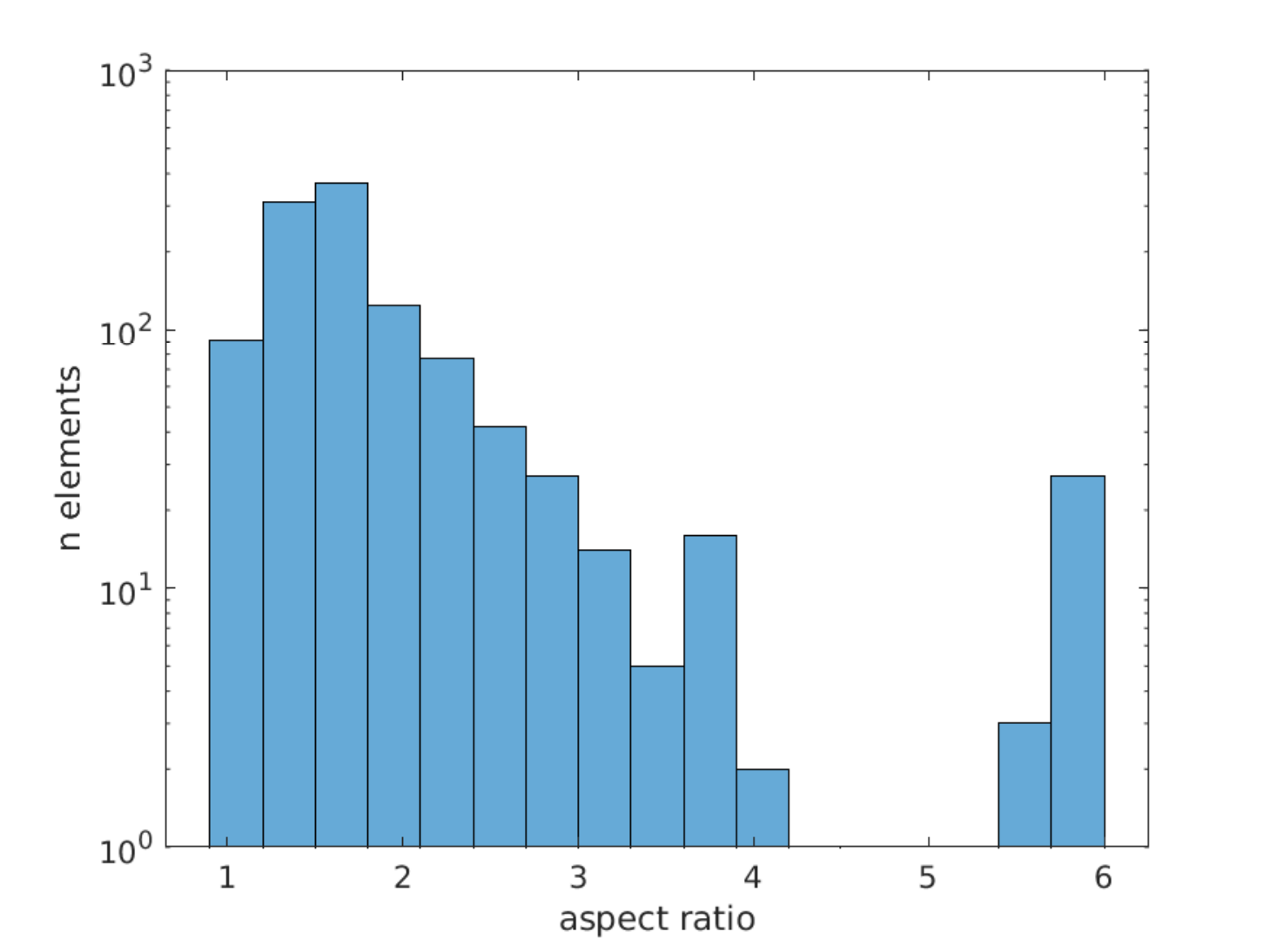}
}
\caption{Histograms of the cells aspect ratio for the different types of grid in test case \ref{subsec:benchmark}.}
\label{fig:anisotropy}
\end{figure}

Another result of the coarsening is a reduction of the number of very
small or very stretched cells. In Figure \ref{fig:anisotropy} we can observe
histograms of an estimate of the cells aspect ratio for the different grids. We
can see that for the Cartesian cut grid and the Voronoi grid the maximum aspect
ratio decreases remarkably with the agglomeration, while in the case of a
Delaunay grid we have the opposite effect. As we will show later {high
anisotropy can result in a less effective stabilization for the MVEM matrix.}
Moreover, in Table \ref{tab:aree_facce} we show that cells agglomeration leads
to an increase of the mean and minimum cell areas, but also to an increase of
the number of faces per cell.

\begin{table}\centering
 \begin{tabular}{|l|c|c|c|c|c|c|}
 \hline
  \multirow{2}{*}&   \multicolumn{3}{c|}{cell area} &\multicolumn{3}{c|}{$n_\mathrm{faces}$}\\     \cline{2-7}
    & average & min & max & average & min & max\\ \hline
    Delaunay & 7.8431e-04 & 8.4186e-05  &2.1020e-03& 3 &  3 & 3\\ \hline
    Delaunay coarse & 9.1075e-04 & 3.9631e-04  & 2.1767e-03 & 3.1557 &  3 & 8\\ \hline
    Cut & 7.7160e-04 &8.4664e-08  &9.1833e-04 &3.9769 &  3 & 6\\ \hline
    Cut coarse& 9.4967e-04 &  3.9945e-04  &2.2589e-03 & 4.4311 &  3 & 10\\ \hline
    Voronoi &  6.5746e-04&  4.6260e-07  &  1.2686e-03 &4.4694 &  3 & 14\\ \hline
    Voronoi coarse &   9.0171e-04&  3.3000e-04 & 3.4502e-03 &5.1109 &  4 & 16\\ \hline
 \end{tabular}
\caption{Average, minimum and maximum value of cell area and number of faces per cell for the six grids employed for test case \ref{subsec:benchmark}}\label{tab:aree_facce}
\end{table}

Referring to the colour code given in Figure \ref{fig:anna_domain}, we set
the aperture $\epsilon = 10^{-4}$ for all the fractures and the permeability
is set to $k_\gamma = \kappa_\gamma = 10^4$ for all the fractures depicted in
red and $k_\gamma = \kappa_\gamma = 10^{-4}$ for the ones in blue. The former
behave as high flow channels while the latter as low permeable barriers. The
rock matrix is characterized by a unit scalar permeability. In
\cite{Flemisch2016a} two sets of boundary conditions were considered,
left-to-right and bottom-to-top. In our case we choose the former, meaning that
we set a value of pressure equal to 4 on the left side of $\Omega$ and to 1 on
the right side of $\Omega$. The
other two boundaries are considered as impervious.

\begin{figure}[htb]
    \centering
    \subfloat[Pressure over line $(0,0.5)-(1,0.9)$\label{fig:anna_pol_pressure}]{
    \includegraphics[width=0.425\textwidth]{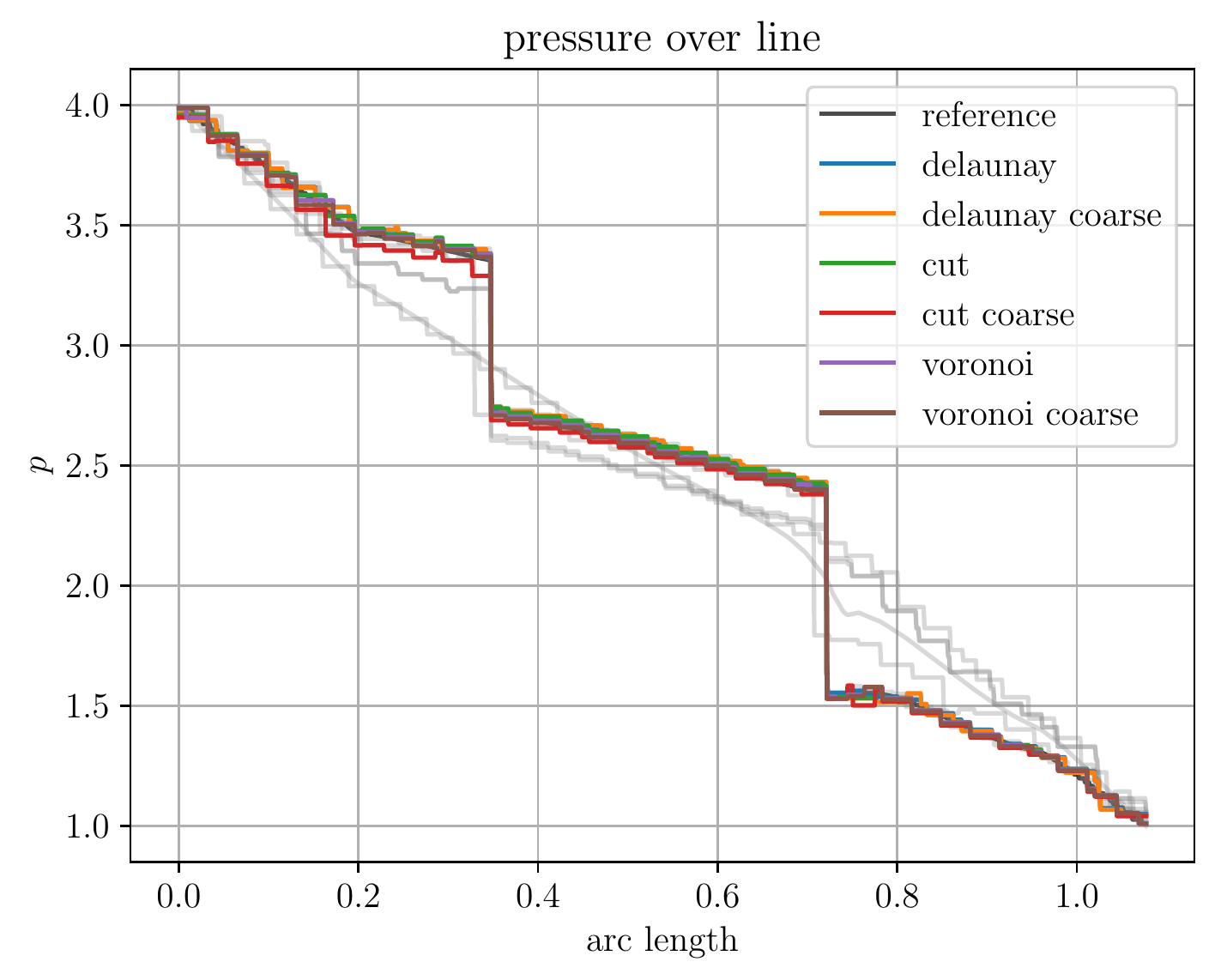}}%
    \hspace{0.1\textwidth}%
    \subfloat[Pressure difference over line $(0,0.5)-(1,0.9)$\label{fig:anna_pol_pressure_diff}]{
    \includegraphics[width=0.425\textwidth]{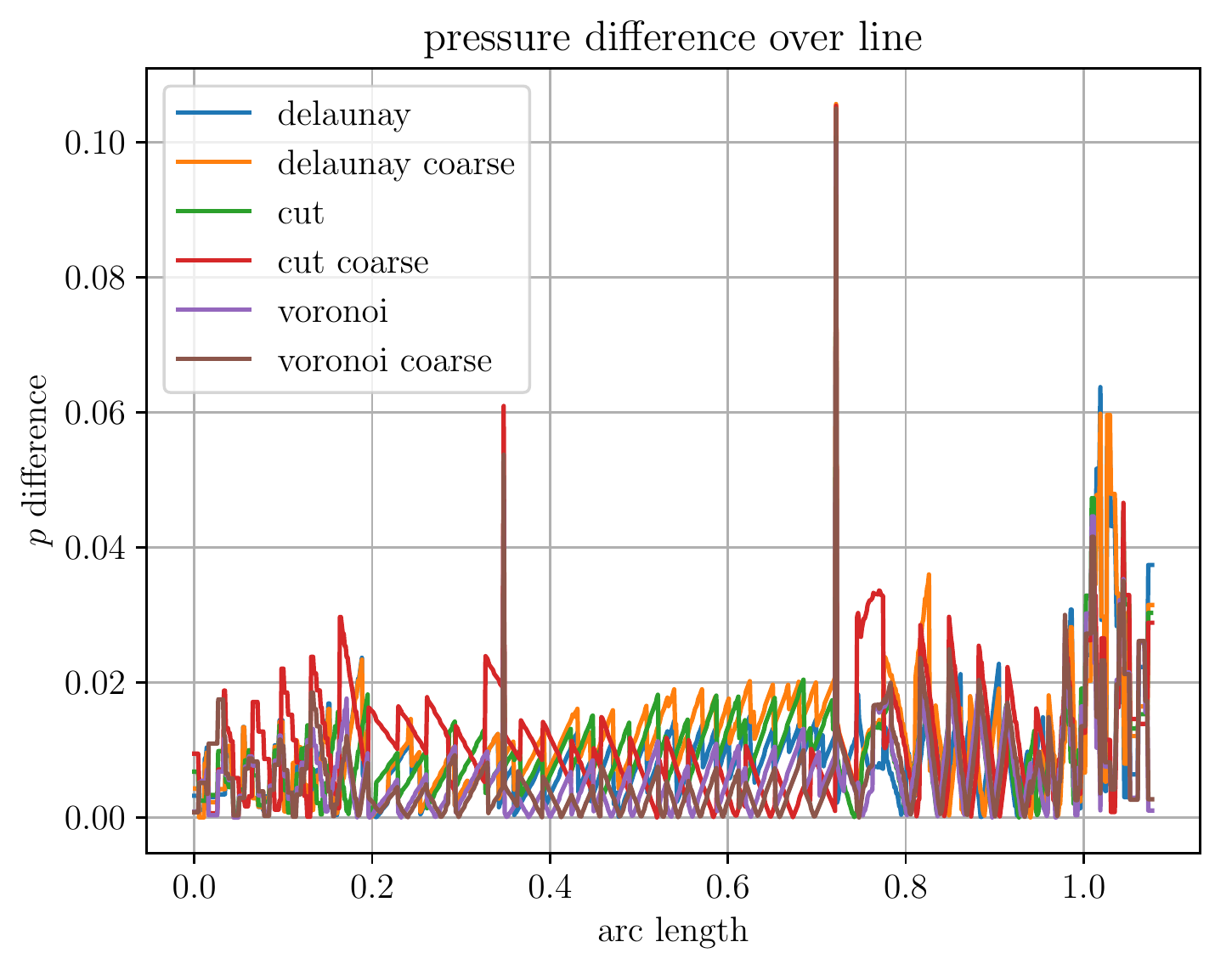}}
    \caption{On the left, pressure over line for the test case of \ref{subsec:benchmark}. The
    grey solutions are the one reported in \cite{Flemisch2016a}. Most of the
    plots overlap with the reference solution, in black. On the right, the
    difference over the same line between a solution and the reference one.}
    \label{fig:anna_pol}
\end{figure}

In Figure \ref{fig:anna_pol_pressure} we report the plot of pressure over the line $(0, 0.5) - (1, 0.9)$, by
using the grids shown in Figure \ref{fig:anna_grids}.
In light grey we present the results obtained in the benchmark
\cite{Flemisch2016a} and in black the reference solution. We clearly see that all
the proposed methods overlap with the reference solution showing high accuracy
even on such coarse grids. In particular results do not deteriorate with the coarsening procedure. Moreover, comparing with the results obtained in
\cite{Flemisch2016a} the ones given by the MVEM are, generally, of higher
quality.

In Figure \ref{fig:anna_pol_pressure_diff}, we show the pressure
difference between the reference solution and the ones obtained with the
considered grids, over the reference solution itself. The errors are quite small
except for the two peaks in
correspondence of the pressure jump of Figure \ref{fig:anna_pol_pressure}. The
reason can be associated to the sampling procedure used in the extraction of
these data.

Finally, as done in \cite{Flemisch2016a} we compute the errors in the rock matrix
between the reference and the computed solution. We consider the following
formula
\begin{gather}\label{eq:error_benchmark}
    err_m^2 = \frac{1}{\abs{\Omega}(\Delta p_{ref})^2}\sum_{f = K_m \cap K_{ref,
    m}} \abs{f} \left(
    p_m|_{K_m} - p_{ref}|_{K_{ref,m}} \right)^2,
\end{gather}
where $p_m|_{K_m}$ is the pressure of the $m$-method at cell $K_m$, $p_{ref}$
is the reference pressure at cell $K_{ref, m}$, and $\Delta p_{ref}$ is the maximum
variation of the pressure on all the domain.
These errors are reported in Table \ref{tab:error_benchmark}.
\begin{table}[htbp]
    \centering
    \begin{tabular}{|l|c|c|}
        \hline
        & original & coarsened\\
        \hline
        Delaunay & 0.013008 & 0.014267 \\
        \hline
        Cartesian cut & 0.012865  & 0.025827\\
        \hline
        Voronoi & 0.0085291 & 0.010037\\
        \hline
    \end{tabular}
    \caption{Pressure error between the reference solution and the compute with the MVEM
    by using formula \ref{eq:error_benchmark}.}
    \label{tab:error_benchmark}
\end{table}
All the errors are quite small and comparable with those reported in
\cite{Flemisch2016a}. When the coarsening procedure is adopted, the errors
slightly increase due to the smaller number of cells except for the Cartesian cut
case where the error doubles, remaining nevertheless acceptable.

Let us now analyze the properties of the system matrix to verify what is the
impact of element size and shape in the different cases. We remind that the
grids have been generated with comparable resolution to obtain similar
numbers of degrees of freedom,  however, the number of unknowns is not exactly
the same. Results are summarized in Table \ref{tab:matrici_caso2}. From the
point of view of the degrees of freedom the Voronoi grid is the most demanding
because, for a given space resolution it generates very small cells close to the
intersections and tips, however, it is also the one that benefits the most from
coarsening. The conditioning is of the same order of magnitude for all grids,
and improves with coarsening. In particular the best result is obtained for the
coarsened Voronoi grid despite the large number of faces per element that
results from clustering of general polygons and reflects in the slightly larger number of non-zero entries per row.

\begin{table}\centering
 \begin{tabular}{|l|c|c|c|c|c|}
 \hline
    & $N_\mathrm{DOF}$ & $N_\mathrm{cells}$ & $N_\mathrm{faces}$ & $\overline{n}$ & $K(A)$ \\ \hline
    Delaunay & 3741 & 1373 &2162 & 5.15 &  4.82E+10 \\ \hline
    Delaunay coarse &3384 & 1196  & 1982 & 5.51 &  3.85E+10 \\ \hline
    Cut & 4961 &1495  &1296 & 6.00 & 4.23E+10 \\ \hline
    Cut coarse& 4474 &  1252 &2814& 7.42 &  3.67E+10 \\ \hline
    Voronoi & 6095 & 1738  & 3913 & 7.33 & 4.10E+10  \\ \hline
    Voronoi coarse & 5118&1326 & 3348 & 9.32&  3.21E+10
 \\ \hline
 \end{tabular}
\caption{Matrix properties for Test case \ref{subsec:benchmark}}\label{tab:matrici_caso2}
\end{table}

\begin{figure}[htb]
    \centering
    \subfloat[Delaunay grid]{
     \includegraphics*[width=0.45\textwidth]{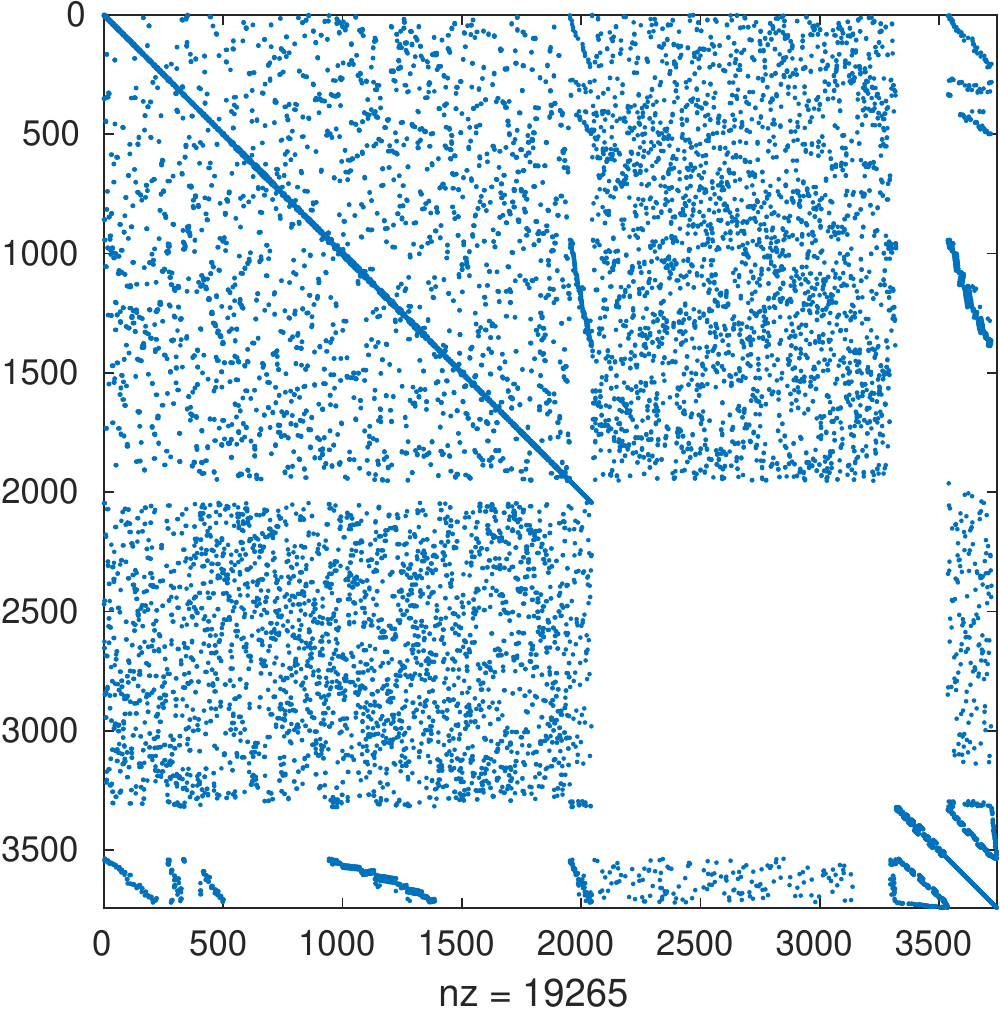}}%
    \hspace{0.05\textwidth}%
    \subfloat[Cartesian cut grid]{
    \includegraphics*[width=0.45\textwidth]{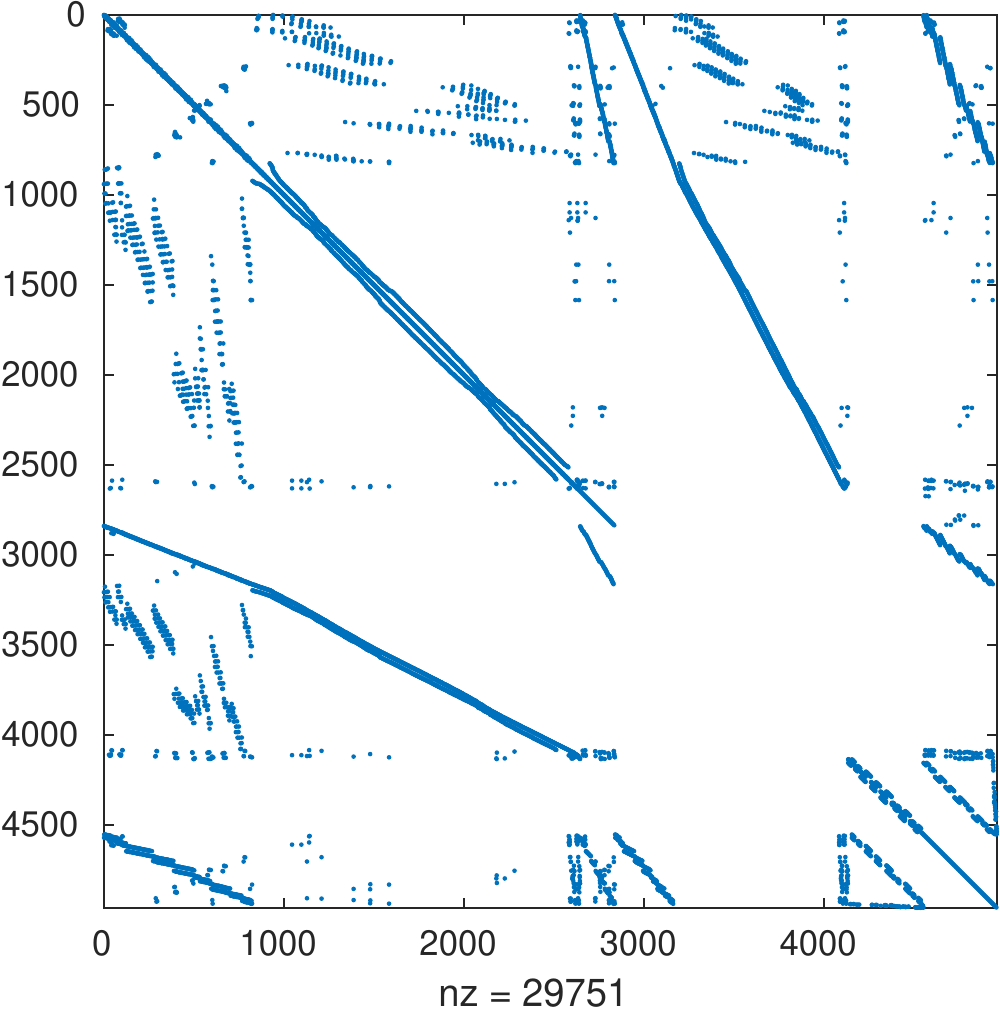}}
    \caption{The sparsity patterns for a Delaunay grid (left) and the Cartesian cut grid (right). }
    \label{fig:sparsity}
\end{figure}

We can also observe that, even if the sparsity of the matrices is similar in all
cases, the pattern can change significantly. In Figure \ref{fig:sparsity} we
compare the matrix structure corresponding to a Delaunay grid and a Cartesian
cut one: the underlying structure of the Cartesian grid has a visible impact on
the sparsity pattern. A similar structure is observed for the case of the
Voronoi grid since, away from the fracture network, the seeds are positioned to
obtain a Cartesian grid. Let $T_\alpha$ denote the time required to solve 1000
times the system arising from the discretization on a mesh $\alpha$ with the
''$\backslash$'' method from \texttt{Matlab}\textsuperscript\textregistered, and
let $\tilde{T}_\alpha=\frac{T_\alpha}{(N_\mathrm{DOF}^\alpha)^3}$ be the time
normalized against the third power of the system size. The corresponding values,
reported in Table \ref{tab:tempi}, seem to indicate that, for the same sparsity,
a faster solution is obtained with a more compact pattern. Solution
strategies for this kind of problem can be found in \cite{Dassi2019a}.

\begin{table}\centering
 \begin{tabular}{|l|c|c|c|c|c|c|}
 \hline
  & Delaunay & Delaunay coarse & Cut & Cut coarse & Voronoi & Voronoi coarse\\ \hline
 $\tilde{T}_\alpha$ &2.630e-10& 3.410e-10   &1.889-e10& 2.262e-10   &1.394e-10& 2.246e-10\\ \hline
 \end{tabular}
 \caption{Normalized time for the solution of the linear systems corresponding to the different grids.}\label{tab:tempi}
 \end{table}

Finally, we study the effect of element shape on the MVEM stabilization term. We define element-wise an index
\begin{gather*}
    \kappa_i=\frac{||S_i||}{||S_i||+||A_i||}
\end{gather*}
where $S_i$ and $A_i$ are the local stability and consistency contributions to the matrix arising from the discretization of the bilinear form $a_\Omega$ on the $i-$th element.

\begin{figure}[htb]
    \centering
    \subfloat[Delaunay grid]{
    \includegraphics[width=0.45\textwidth]{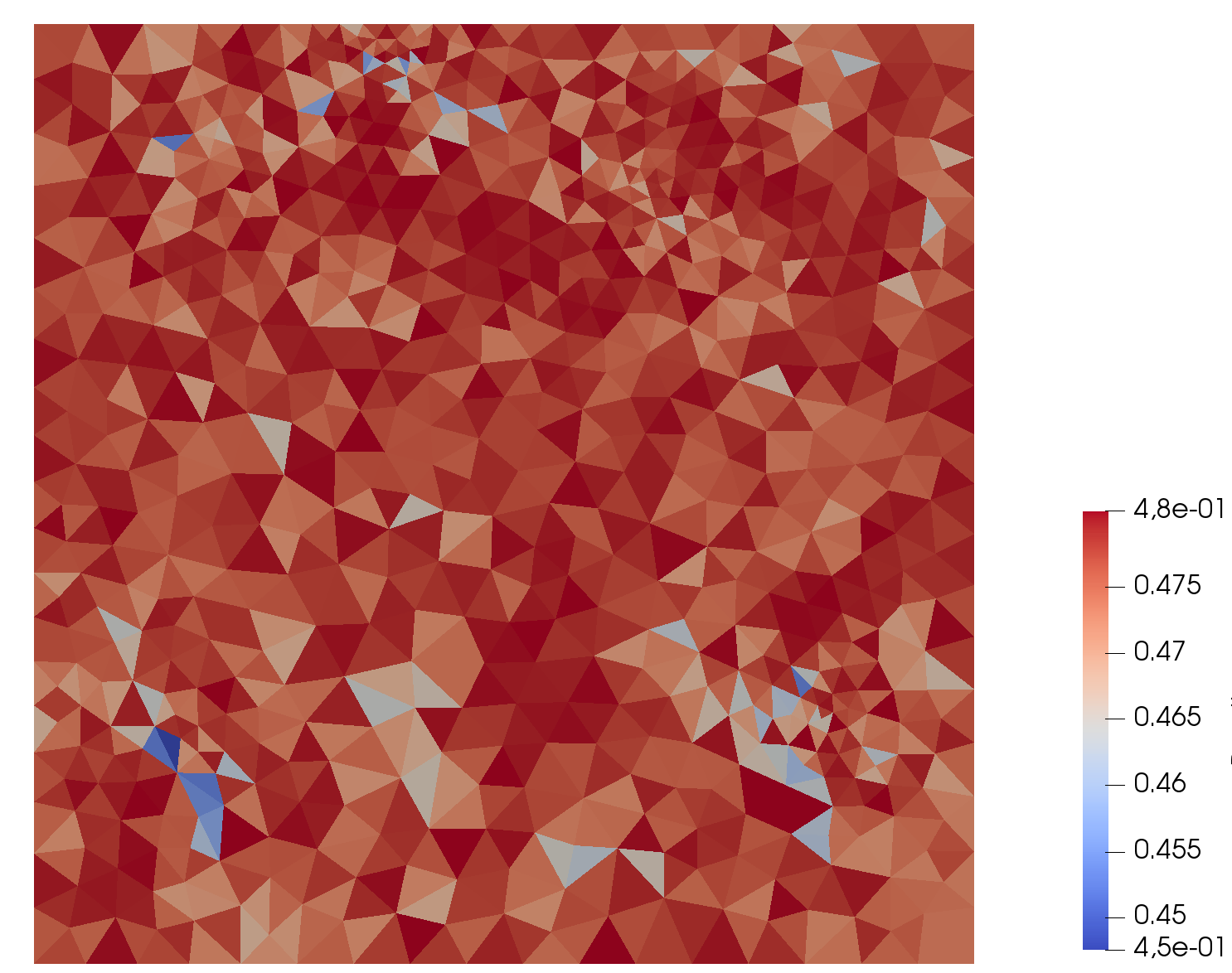}}%
    \hspace{0.05\textwidth}%
    \subfloat[Voronoi grid]{
    \includegraphics[width=0.45\textwidth]{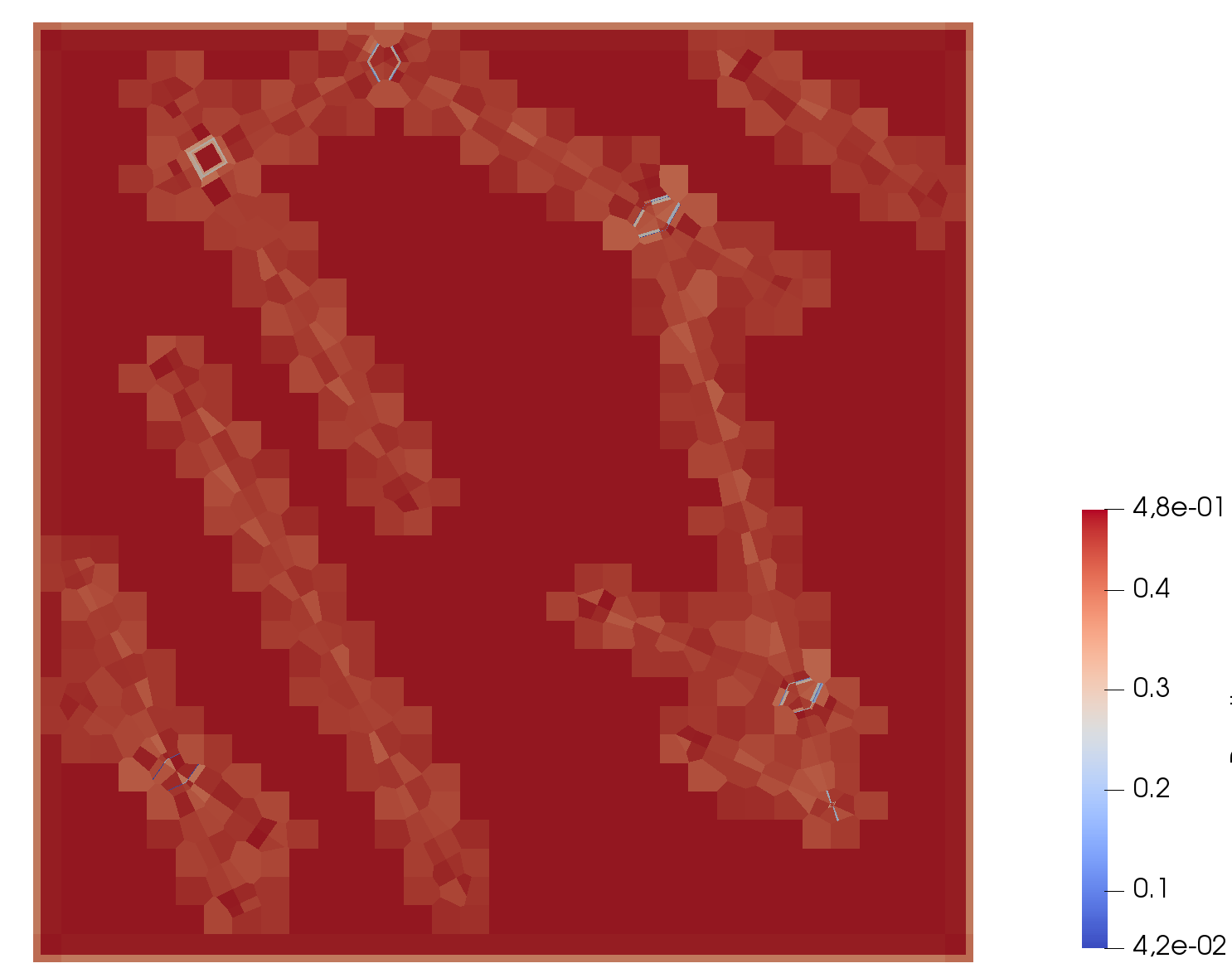}}
    \caption{On the left, $\kappa_i$ on the Delaunay grid, on the right the same
    index on the Voronoi grid before clustering. }
    \label{fig:stab}
\end{figure}

As shown in Figure \ref{fig:stab} in the case of a Delaunay grid the norm of the
stabilization term in each local matrix is comparable to the norm of the
consistency term, i.e. $\kappa_i\simeq 0.5$ everywhere. In the Voronoi grid
instead we have elements with extremely high aspect ratios (up to 60), or, in
other words, we have small edges compared to the typical mesh size. In this
latter case the norm of the stability term is one order of magnitude smaller in
elements with very small edges. A discussion of the stability bounds for grids in
the case of small edges can be found in \cite{BeiraoVeiga2016},
\cite{BeiraoVeiga2017} for the primal formulation of elliptic problems.

\section{Conclusion}\label{sec:conclusion}

In this work we have presented and discussed the performances of the Mixed
Virtual Finite Element Method applied to underground problems. One of its main
advantages is the possibility to handle, in a natural way, grid cells of any
shape becoming suitable for its usage in problems with complex geometries, such
as subsurface flows. A second strong point is the ability of the scheme to
handle, in a robust way, strong variations of the permeability matrix which is
again a common aspect for underground processes.  Finally, the numerical scheme
is also locally mass conservative making it very suitable in the coupling of
other physical processes, like transport problems.  We have tested the
capabilities of the scheme with respect to two test cases that are known in
literature and stress the two aforementioned critical points: heterogeneity and
geometrical complexity. A first remark is that the mixed virtual element method
gives high quality results also for challenging grids and physical data, making
it a promising and interesting scheme for industrial applications. Moreover we
performed some comparisons of the system matrices arising from the
discretization of the problem on different types of grids: Delaunay, Voronoi,
Cartesian grids cut by fractures. We observed similar condition numbers and
sparsity, but a better sparsity patterns for grids obtained from the
modification of structured ones. We also applied coarsening by means of
permeability based and volume based clustering: besides reducing the
computational cost this technique allowed us to eliminate small cells and, in
some cases, cells with very large aspect ratios where the MVEM stabilizazion term
employed in this work does not scale correctly. {Future research may focus on
the choice of the most stabilization term formulation for the grid type, as well
as to the generalization of this work to the three dimensional case, including
the discussion of corner point grid which are widely used in subsurface flows
but pose many challenges due to the presence of non-planar faces and non-convex
elements.}

\section{Acknowledgements}

We acknowledge the PorePy development team:
Eirik Keilegavlen, Runar Berge, Michele Starnoni, Ivar Stefansson, Jhabriel Varela,
Inga Berre.


\end{document}